\documentclass{article}
\usepackage{graphicx} 
\usepackage{algorithm}
\usepackage{algpseudocode}
\usepackage{todonotes,pdfpages}

\newcommand{\R}{\mathbb{R}}	
\newcommand{\SDP}{\mathcal{S}}

\def\trace{\mbox{trace }}

\def\ddd{{\mathsf d}}
\def\e{{\mathsf e}}

\def\ooo{{\mathsf o}}

\def\uu{{\mathsf u}}
\def\vv{{\mathsf v}}
\def\w{{\mathsf w}}
\def\x{{\mathsf x}}

\def\Rb{{\mathsf R}}

\def\Ub{{\mathsf U}}
\def\Vb{{\mathsf V}}
\def\Wb{{\mathsf W}}
\def\Xb{{\mathsf X}}

\def\Zb{{\mathsf Z}}
\def\irg#1#2{{ \in\! [{#1}\! : \! {#2}]}}
\newcommand{{\T}}{{^\top}}

\def\diag{\mbox{diag }}

\def\beq#1{\begin{equation}\label{#1}}
\def\eeq{\end{equation}}
\def\bep{\begin{proof}}
\def\ep{\end{proof}}

\def\R{{\mathbb R}}

\def\sss{{\mathsf s}}

\parskip=\medskipamount
\title{Feature selection in linear Support Vector Machines\\ via a hard cardinality
 constraint:\\ 
 a scalable conic decomposition approach}
\author{Immanuel Bomze \thanks{Faculty of Mathematics and Research Network Data Science, University of Vienna, Oskar-Morgenstern-Platz 1,
1090 Wien, Austria. ORCID ID: 0000-0002-6288-9226 E-mail: {\tt immanuel.bomze@univie.ac.at}}\and Federico D'Onofrio \thanks{Department of Computer, Control and Management Engineering, Sapienza University of Rome, Via Ariosto 25, 00185 Rome, Italy. ORCID ID: 0000-0002-0584-5534 E-mail: {\tt federico.donofrio@uniroma1.it}} \and Laura Palagi \thanks{Department of Computer, Control and Management Engineering, Sapienza University of Rome, Via Ariosto 25, 00185 Rome, Italy. ORCID ID: 0000-0002-9496-6097 E-mail: {\tt laura.palagi@uniroma1.it}} \and Bo Peng\thanks{VGSCO, University of Vienna, Oskar-Morgenstern-Platz 1,
1090 Wien, Austria. ORCID ID: 0000-0002-2650-0295 E-mail: {\tt bo.peng@univie.ac.at}}}
\date{\today}

\usepackage{hyperref}
\usepackage{amsfonts}
\usepackage{geometry}

\usepackage{rotating}
\usepackage{amssymb}   
\usepackage{amsthm}
\usepackage{amsmath}

\usepackage{pict2e}
\usepackage{keyval}
\usepackage{calc}
\usepackage{fp}
\usepackage{diagbox}
\usepackage{booktabs}
\usepackage{tabularx,colortbl, makecell, caption, multirow}
\usepackage{boldline}

\newcommand{\red}{\textcolor{black}}

\def\T{\mathcal{T}}

\newcommand{\Z}{\mathbb{Z}}
\def\sgn{\mbox{\rm sgn}}

\newtheorem{theorem}{Theorem}[section]
\newtheorem{proposition}[theorem]{Proposition}
\newtheorem{lemma}[theorem]{Lemma}

\newtheorem{corollary}[theorem]{Corollary}

\newtheorem{remark}[theorem]{Remark}

\usepackage{rotating}
\usepackage{graphicx}
\usepackage{bm}
\usepackage{xcolor}
\usepackage{multirow}
\usepackage{nicematrix}

\usepackage{graphicx}
\graphicspath{{images/}}
\usepackage{float}
\usepackage{longtable}
\usepackage{adjustbox}

\newcommand{\ds}{\displaystyle}

\usepackage{colortbl} 

\graphicspath{{./figures/}}

\usepackage[titletoc,toc,title]{appendix}
\usepackage{comment}
\usepackage{graphics}
\usepackage{graphicx}
\usepackage{adjustbox}
\usepackage{tikz}
\usepackage{caption}
\usepackage{psfrag}
\usepackage{soul}
\usepackage{mwe}
\usepackage{pgfplots}
\pgfplotsset{compat=1.17} 

\graphicspath{{./figures/}}

\usepackage{comment}
\usepackage{graphics}
\usepackage{graphicx}
\usepackage{adjustbox}
\usepackage{tikz}
\usepgfplotslibrary{dateplot}
\usepgfplotslibrary{groupplots}
\usepackage{caption}
\usepackage{subcaption}
\usepackage{psfrag}
\usepackage{soul}
\usepackage{mwe}
\usepgfplotslibrary{fillbetween}

\usetikzlibrary{er,positioning, calc, shapes, arrows, chains, fit}

\definecolor{mygreen}{RGB}{80,176,50}

\usepackage{xcolor}

\makeatletter

\newcommand*{\@rowstyle}{}

\newcommand*{\rowstyle}[1]{
  \gdef\@rowstyle{#1}%
  \@rowstyle\ignorespaces%
}

\newcolumntype{=}{
  >{\gdef\@rowstyle{}}%
}

\newcolumntype{+}{
  >{\@rowstyle}%
}

\makeatother

\definecolor{dred}{rgb}{0.55, 0.0, 0.0}
\definecolor{dblue}{rgb}{0.0, 0.2, 0.6}
\definecolor{dgreen}{RGB}{80,176,50}

\usepackage{colortbl}

\definecolor{mygreen}{RGB}{80,176,50}
\definecolor{myred}{RGB}{222,36,16}
\definecolor{myind}{RGB}{93,174,255}
\definecolor{mypink}{RGB}{219,129,255}
\definecolor{azzdark}{RGB}{42,72,180}
\definecolor{azzlight}{RGB}{53,146,203}
\definecolor{verdone}{RGB}{0,128,0}
\definecolor{indaco}{RGB}{68,137,206}
\definecolor{celestino}{RGB}{222,247,254}
\definecolor{verdeacceso}{RGB}{71,185,109}
\definecolor{pred}{RGB}{182,7,40}
\definecolor{pblue}{RGB}{58,76,189}

\usetikzlibrary{er,positioning, calc, shapes, arrows, chains, fit, decorations.pathreplacing}

\definecolor{orange}{rgb}{0.8, 0.47, 0.13}
\definecolor{green}{rgb}{0.0, 0.62, 0.42}

\begin{document}

\maketitle
\begin{abstract}

  In this paper, we study the embedded feature selection problem in linear Support Vector Machines (SVMs), in which a cardinality constraint is employed,  leading to \red{an interpretable classification model}. The problem is NP-hard due to the presence of the cardinality constraint, even though the original linear SVM amounts to a problem solvable in polynomial time. To handle the hard problem, we first introduce two mixed-integer formulations for which novel semidefinite relaxations are proposed. Exploiting the sparsity pattern of the relaxations, we decompose the problems and obtain equivalent relaxations in a much smaller cone, making the conic approaches scalable. To make the best usage of the decomposed relaxations, we propose heuristics using the information of its optimal solution. Moreover, an exact procedure is proposed by solving a sequence of mixed-integer decomposed semidefinite optimization problems. Numerical results on classical benchmarking datasets are reported, showing the efficiency and effectiveness of our approach. 
\end{abstract}

\section{Introduction}

In supervised classification, Machine Learning (ML) models are trained on labeled datasets to learn how to predict the correct output for new, unseen data. Classification algorithms are prone to overfitting when the dataset is small and contains irrelevant features, or noise~\cite{hastie2009elements}.  This problem can also occur with correlated features or features containing errors.
 Further, interpretability in ML is gaining more importance, and remains a cutting-edge research topic. Indeed, designing models which provide explanatory insights and help end users to understand the rationale behind decisions, is crucial from both predictive and debugging perspective~\cite{kamath2021explainable}.

It is therefore beneficial to use a smaller set of features that can give a good approximation of the target function, i.e., the classifier. Feature Selection (FS) aims to eliminate less significant and/or redundant features to identify the possibly optimal subset for the classifier, which becomes crucial when dealing with real data.  Among various feature selection techniques, embedded methods integrate feature selection within the classifier training process.
In this paper, we study the embedded feature selection problem for SVMs.

\red{From now on, we use lowercase letters either in the sans-serif Latin font, e.g., $\mathsf{x}$, or in the boldface font for Greek letters, e.g., $\bm{\xi}$, to indicate vectors.}
\red{Transposition is denoted by $\hbox{}^\top$ and the zero vector  by $\ooo$; given two integers $a,b$ with $a\le b$,  we abbreviate by $[a\! : \! b]=\{k\in \Z:\ a\leq k \leq b\}$.}

Given a dataset with $m$ samples, and $n$ features, $\{(\x^i,y^i) \in \mathbb{R}^n\times\{-1,1\}: i\irg 1m\}$,  the linear SVM  classification problem defines a classifier as the \red{sign} function 
$f:\R^n\to \{-1,1\}$, 
$$f(\x) = \sgn({\w^*}^\top \x + b^*)\, ,$$
  using the structural risk minimization
principle~\cite{vapnik1999nature}.

In traditional SVMs~\cite{VapnikSVM1995}, coefficients $(\w^*,b^*)\in\R^n\times\R$ identify a separating hyperplane (if it exists) that maximizes its margin, i.e. the distance between the hyperplane and the closest data points of each class. Accounting for classes not linearly separable, the tuple  $(\w^*,b^*)$ is found by solving the following convex quadratic problem:

\begin{equation}\label{prob:l2-SVM} 
\begin{array}{rl}
    (\text{SVM}) \quad
    \min \limits_{{\w,b,\bm{\xi}}} \quad & \ds\frac 1 2  \|\w\|_2^2 +C \sum_{i=1}^m \xi_{i} \\
    \text{s.t.} \quad 
    & y^i(\w^\top \x^i+b)\ge  1 - \xi_i\, ,\,   \quad i\irg{1}{m}\, ,  \\
    & \bm{\xi} \geq  \ooo\, .
\end{array}
\end{equation}

The non-negative variables $\xi_i$  allow points to be on the wrong side of their “soft margin”, as well as being "too close" to the decision boundary.

The sum $\sum_{i = 1}^m\xi_i$ in the objective function represents an upper bound on the total misclassification error, which is weighted by the hyperparameter $C\!>\!0$ so to balance the maximization of the margin, which is inversely proportional to $\|\w\|_2$.  Problem~\eqref{prob:l2-SVM} is categorized as $\ell_2$-regularized $\ell_1$-loss SVM  (see~\cite{fan2008liblinear,hastie2004entire}, for instance). When needed, we can also refer to it as $\ell_2$-SVM to highlight the use of the $\ell_2$ regularization contrasting with the use of $\ell_1$ regularization, as in~\cite{hastie2004entire,nguyen2010optimal}, for instance.

The optimal solution $\w^*$ of \eqref{prob:l2-SVM} is usually dense, in the sense that usually all components of $\w^*$ are nonzero and thus every feature contributes to the definition of the classifier. Indeed, through duality theory, $\w^*$ can be expressed as
the linear combination of a subset of samples, called support vectors. More precisely, $\w^*=\sum_{i = 1}^m \alpha^*_i y^i \x^i$, with $\alpha_i^* \ge 0$ being the dual optimal solution of problem~\eqref{prob:l2-SVM} with nonzero value only if sample $i$ represents a support vector. Thus, it is unlikely to obtain $w_j^*=0$ for some feature $j$, and sparsity in the vector $\w^*$ must be forced.

In this paper, we aim to strictly control sparsity by imposing hard cardinality constraints on the number of nonzero components in $\w^*$. We approach the feature selection problem from a combinatorial perspective, using binary variables to clearly indicate which features are selected.
More specifically, the SVM problem with feature selection we are interested in is: 
\begin{equation}
\label{prob:l2-FS-SVM} 
\begin{array}{rl}
    (\text{FS-SVM}) \quad
    \min \limits_{{\w,b,\bm\xi} } \quad &  \ds\frac 1 2 \| \w\|^2_2 + C \sum_{ i=1}^m \xi_{i}\\ 
       \text{s.t.} \quad
    & y^i ( \ds \w^\top \x^i + b) \ge 1 - \xi_{i} \, ,  \quad i\irg{1}{m} \, ,  \\
    &\|\w\|_0\leq B  \\
    & \bm\xi \geq \ooo\, , \\
\end{array}
\end{equation}
where the $\ell_0$-pseudonorm is defined as $\| \w\|_0 := \sum_{j=1}^{n} |\sgn(w_j)|$. 
The user-defined parameter $B\in\mathbb{Z}_+$ constitutes a budget on the number of features used, indicating that the number of features involved in the classifier is at most $B$. Even though the objective function and all other constraints are convex quadratic or linear, the cardinality constraint $\|\w\|_0\leq B$ is of combinatorial type, which renders~\eqref{prob:l2-FS-SVM} NP-hard. \red{The class of optimization problems involving cardinality constraint has gained significant attention over the past decades, for a comprehensive overview of such problems, we refer the interested readers to \cite{bertsimas2009algorithm, gao2013optimal, kanzow2021augmented, tillmann2024cardinality} and references therein.}

\red{By strictly controlling the budget value 
$B$, we can achieve a sparse, more transparent model which is usually easier to explain. Sparsity enhances interpretability by limiting the number of features, which simplifies the model structure and clarifies its decision-making, allowing clearer insights into predictions \cite{rudin2022interpretable}. A model that selects at most
$B$ features offer precise control over sparsity, ensuring alignment with user-defined complexity requirements. Unlike penalty-based methods, such as Lasso regularization, which incorporates a $\ell_1$-norm penalty in the objective function, this approach provides a known feature count {\em a priori}, enhancing computational efficiency and reducing the need for extensive tuning \cite{bertsimas2020sparse}. Additionally, models with cardinality constraints are often generalized to budget-constrained models, where features are selected based on both relevance and cost, creating solutions that are sparse and cost-effective \cite{labbe2019mixed, lee2020mixed}.} 

\red{Although the performance of linear SVMs may fall short in certain classification tasks when compared to more complex models, the use of feature maps into high-dimensional feature spaces somehow conflicts with the request of interpretability nowadays considered as a crucial aspect in machine learning; indeed, the original variables have a direct meaning more easily understood by practitioners, particularly if only a few of them (regulated by $B$) are used in the classifier.}

The main contributions of our work are:
\begin{itemize}
    \item[--] We \red{introduce} two novel Mixed Integer Quadratic optimization Problem (MIQP) formulations for the FS-SVM problem \eqref{prob:l2-FS-SVM},  \red{building upon on known optimization techniques for   
    tackling the $\ell_0$-pseudonorm constraint,  either by means of a big-M reformulation or by modelling it with complementarity constraints}.
    \item[--] We present novel relaxations that decompose with respect to small blocks of variables, thus being easily solvable when dealing with datasets with a large number of features.
    \item[--] We propose both exact and heuristic algorithms that are easily implementable and exploit the proposed relaxations.
    \item[--] We conduct a battery of extensive numerical experiments that demonstrate the effectiveness of our algorithms in solving the proposed problem compared to off-the-shelf solvers, \red{and validate the model from a machine learning perspective}.
\end{itemize}

\subsection {Related work}

Embedding feature selection in linear SVM has been proposed in several papers. One of the most frequently used approaches consists in using the $\ell_1$-regularized SVM, the  $\ell_1$-SVM, instead of the $\ell_2$-SVM. Indeed, the $\ell_1$-regularization, also known as LASSO regularization, leads to sparser solution and to \red{an} easier problem, being formulated as a Linear Programming problem
(see~\cite{FungFSNewton,hastie2004entire,nguyen2010optimal,yuan2010comparison}, for instance).

While replacing the $\ell_2$-regularization by the $\ell_1$-regularization has a significant impact on reducing the computational requirement for solving the optimization problem, the $\ell_1$-SVM leads to relatively reduced classification accuracy due to the inability to maximize the margin of separation between the two classes of data.  To overcome the disadvantages, 
in~\cite{ wang2006doubly,zou2005regularization} the elastic-net regularization technique has been proposed which uses both the $\ell_2$- and the $\ell_1$-regularization in the objective function. 

\red{Another approach penalizes the $\ell_0$-pseudonorm of $\w$ in the objective function, offering a better approximation than the $\ell_1$-norm.}
In this framework, 
\cite{Bradley1998FeatureSV} studied the model in which the $\ell_2$-regularization is substituted with the $\ell_0$-pseudonorm of $\w$. They solve the resulting problem by approximating $\|\w\|_0 \approx 1-e^{-\lambda \|\w\|_1}$, where $\lambda \geq 0$. The resulting problem is the minimization of a concave function over a polyhedral set and the authors proposed a successive linearization algorithm 
to find a local solution.  \cite{weston2003use} tackled the the same problem as~\cite{Bradley1998FeatureSV} and proposed to solve it by iteratively rescaling the input features values, multiplying them by the absolute value of the weights $w_j$, obtained from the resolution of the SVM problem, until convergence is reached. Variables can be ranked by removing those features whose weights become zero during the iterative algorithm and computing the order of removal. 

None of \red{the} approaches described above have a direct and hard control on the number  $\|\w\|_0$ of nonzero components in the vector $\w^*$ like the property  $\|\w\|_0\le B$.

Facing the problem~(FS-SVM), \cite{chan2007direct} proposed to replace the cardinality constraint by the non-linear constraint
$\|\w\|_1^2\le B\|\w\|_2^2$. This non-convex constraint can be relaxed, leading to two convex relaxations of the sparse SVM.
However, the resulting optimal solutions seem to be not as strong~\cite{ghaddar2018high}, and they violate the imposed desired features limit $B$. To overcome these drawbacks, \cite{ghaddar2018high} proposed an alternative relaxation, adjusting a parameter value iteratively to obtain improving bounds and employing a bisection algorithm to ensure that $\|\w\|_0$ does not exceed the maximum limit $B$.

Another approach based on bilevel optimization is proposed by \cite{agor2019feature}, where binary variables in the upper-level problem are used \red{to explicitly} control the number of features selected. They proposed a genetic algorithm for the solution of the resulting problem.

 \cite{maldonado2014feature} proposed a Mixed Integer Linear Optimization Problem (MILP) model, introducing feature selection through a budget constraint to limit the number of features used in the classification process, but focusing just on minimizing the sum of the deviations $\xi_i$ rather than considering also the margin. Extending this work, \cite{labbe2019mixed} proposed a MILP formulation which considers also the margin maximization, to select suitable features for constructing separating hyperplanes, incorporating the same budget constraint. The authors proposed both heuristic and exact approaches for solving their formulation when applied to large datasets.

As regards the $\ell_2$-SVM,  \cite{aytug2015feature} proposed two feature selection (hard and soft constrained) models where binary variables act as a multiplicative "mask" on the $\w$, and enter both  the objective function and the constraints in a nonlinear way. An exact method based on Generalized Benders Decomposition is proposed to solve the two formulations. \red{However, the experimental results indicate that the algorithm fails to converge in practice for datasets with a large number of features, reaching the preset limit on the decomposition iterations}. 
Another approach employing $\ell_2$-SVMs and cardinality constraints is used in \cite{MARGOT}, where the authors proposed a new MIQP model for building optimal classification trees with SVM sparse splits. 

\subsection{Paper organization and notation}
This paper is structured as follows: in Section~\ref{Sec2}, two mixed-integer reformulations of the targeted problem~\eqref{prob:l2-FS-SVM} are proposed. Several tractable relaxations of these two mixed-integer problems are introduced in Section~\ref{Sec3} with a theoretical comparison among them. In the same section, we further propose novel decomposed relaxations, exploiting the sparsity pattern of the problems. Section~\ref{Sec4} focuses on algorithmic methods: a heuristic and an exact algorithm are proposed for solving the problem \eqref{prob:l2-FS-SVM}, based on the optimal solutions of the relaxations. Applying these algorithms to some benchmark datasets, 
we report the numerical results in Section~\ref{Sec5}, from the optimization point of view and from the ML perspective. The concluding Section~\ref{Sec6} also discusses possible future research directions.

Throughout the paper, we denote by \red{lower case  letters, e.g.~$\x$, and by upper case letters, e.g.~$\Xb$, a vector and a matrix, respectively.} Especially, $\e$ and $\ooo$ represent all-ones and all-zeros vectors of proper dimension. The dimension should be clear from the context unless otherwise mentioned. We employ subscripts to indicate a specific element of vectors or matrices. For instance, by $x_i$ we denote the $i$-th element of the vector $\x$, and by $X_{ij}$ we denote the $ij$-th entry of the matrix $\Xb$.
Furthermore, by $\diag\Xb=[X_{ii}]_{i=1}^n$ we denote the vector of diagonal of $\Xb$. 
Finally, by $\SDP^n$ and $\SDP_+^n$ we denote the set of  real symmetric $n\times n$ matrices, and the set of symmetric positive-semidefinite matrices of order $n$, respectively. 

Additionally, let us introduce the set
\begin{equation}\label{eq:xi}
\Xi=\{(\w,b,\bm{\xi})\in \R^{n+1+m}: y^i ( \ds \w^\top \x^i + b) \ge 1 - \xi_{i}\, ,  \  i\irg{1}{m}\, , \bm{\xi} \geq \ooo\} \, ,  
\end{equation}
so that problem~\eqref{prob:l2-FS-SVM} reads as:
$$\begin{array}{rll}
    \min \limits_{{\w,b,\bm{\xi}} } \quad &  \ds\frac 1 2 \|\w\|^2_2 + C \sum_{ i=1}^m \xi_{i} \\ 
       \text{s.t.} \quad  & ({\w,b,\bm{\xi} }) \in \Xi  \\[.8em] 
       &\|\w\|_0\leq B.
\end{array}$$

In the next sections, we will also make use of big-M reformulations, so it is convenient to introduce the two sets
\begin{equation}\label{eq:omega}
\Omega_M=\{(\w,\vv)\in \R^{2n}:    -M \vv \leq \w \leq M \vv \,, \, \e^\top \vv = B  \}     
\end{equation}
and 
\begin{equation}\label{eq:Bset}
{\cal B}_M=\{(\w, \uu )\in \R^{2n}: -M(1-\uu)\leq \w \leq M(1-\uu)\, ,\,  \e^\top \uu = n - B \}\, ,    
\end{equation}
where the big-M parameter $M>0$ is large enough
(see Remark~\ref{rem3.2} below) and the variables $\uu=1-\vv$ are the negated $\vv$, sometimes more convenient in notation.

\section{Exact Mixed-Integer approaches}\label{Sec2}

In this section, we propose two exact mixed-integer formulations of the original problem~\eqref{prob:l2-FS-SVM}, using two different mixed-integer systems to present the cardinality constraint with \red{auxiliary} binary variables introduced as \red{indicators}, namely a big-M strategy and a complementarity constraint strategy. 

Differently from \cite{labbe2019mixed}, we study the following MIQP model incorporating the standard $\ell_2$-norm on the hyperplane weights as the regularization term:
\begin{equation}
\label{prob:l2-FS-SVM1}
\begin{array}{rll}
    (\text{BigMP}) \quad
    \min \limits_{{\w,b,\bm{\xi} ,\vv} } \quad &  \ds\frac 1 2 \|\w\|^2_2 + C \sum_{ i=1}^m \xi_{i} \\ 
       \text{s.t.} \quad  & ({\w,b,\bm{\xi} }) \in \Xi  \\[.8em] 
        & ({\w, \vv}) \in  \Omega_M \\[.8em] 
            & \vv \in\{0,1\}^n\, ,
\end{array}
\end{equation}
where $\Xi$ and $\Omega_M$ are defined in~\eqref{eq:xi} and~\eqref{eq:omega}.

To show the equivalence between~\eqref{prob:l2-FS-SVM} and~\eqref{prob:l2-FS-SVM1}, we consider a $B$-sparse $\w$, e.g., $\|\w\|_0 \leq B$. If $\|\w\|_0 = B$, there exists a $\vv\in\{0,1\}^n$ such that $(\w,\vv)\in \Omega_M$, where $\vv$ is given by the rule that $v_j = 0$ if $w_j = 0$ and $v_j = 1$ if $w_j \neq 0$; If $\|\w\|_0 < B$, there also exists a $\vv\in\{0,1\}^n$ such that $(\w,\vv)\in \Omega_M$, where $\vv$ is given by the rule that $v_j = 1$ for $j\in I$ and $v_j = 0$ otherwise, where  $I$ is a index set covers all of features with $w_j \neq 0$, satisfying $|I| = B$ .

\begin{remark}\label{remark1}
    It is worth noticing that the optimal value of model~\eqref{prob:l2-FS-SVM1} remains the same if we substitute the equality constraint $\e^\top \vv = B$ by the budget constraint $\e^\top \vv \leq B$, which is derived directly from the targeted problem~\eqref{prob:l2-FS-SVM}.
\end{remark}

Alternatively, instead of expressing the cardinality constraint by the big-M formulation, we employ a set of binary variables $\uu\in \{0,1\}^n$ and a complementarity constraint to count the nonzero elements of $\w$ and end up with another mixed-integer reformulation of \eqref{prob:l2-FS-SVM}:

\begin{equation}
\label{prob:l2-FS-SVM2}
\begin{array}{rlll}
    (\text{CoP}) \quad
    \min \limits_{{\w,b,\bm{\xi },\uu} } \quad &  \ds\frac 1 2 \|\w\|^2_2 + C \sum_{ i=1}^m \xi_{i} \\ 
       \text{s.t.} \quad
   & ({\w,b,\bm{\xi }}) \in \Xi  \\
    & u_jw_j = 0  \, , &  j\irg{1}{n}\\
    & \e^\top \uu = n - B \\
    & \uu \in\{0,1\}^n  \, .
\end{array}\end{equation} 

 Unlike~\eqref{prob:l2-FS-SVM1}, the binary variables $u_j$ count the number of zero elements of $\w$. For a \eqref{prob:l2-FS-SVM2}-feasible point $(\w,b,\bm{\xi},\uu)$, we let $\vv = 1- \uu$ and the point $(\w,b,\bm{\xi},\vv)$ is \eqref{prob:l2-FS-SVM1}-feasible. To show this, we only need to show $(\w,\vv) \in \Omega_M$, namely, that $(\w,\vv)$ satisfies the big-M constraint in \eqref{eq:omega}. It is clear that the constraint holds for the coordinates of $\w$ with $w_j = 0$. On the other hand,  the complementarity constraint $(1-v_j) w_j = 0$ implies $v_j = 1$ if $w_j \neq 0$, and therefore the point $(\w,\vv) \in \Omega_M$. Hence, the problem \eqref{prob:l2-FS-SVM2} is equivalent to the problem \eqref{prob:l2-FS-SVM1}.
 
\begin{remark}\label{remark2}
     In the above problem \eqref{prob:l2-FS-SVM2}, the constraint $ \e^\top \uu = n - B$ can be equivalently replaced by the constraint $ \e^\top \uu \geq n - B$ for the same reason as Remark~\ref{remark1}, and the binary constraint $u_j\in\{0,1\}$ can be equivalently reduced to the box constraint $0 \leq u_j \leq 1$ due to the existence of the complementarity constraint. To show this, we consider the model with the box constraint and assume that $\w$ is $\rho$-sparse with $\rho > B$. Due to the constraint $u_{j} w_{j} = 0$, $u_j = 0$ if $w_{j}$ is nonzero and $0 \leq u_j \leq 1$ if $w_{j}$ is zero. Hence, $\e^\top \uu \leq n-\rho < n- B$, which contradicts the budget constraint. Therefore all  $\w$ vectors with cardinality strictly larger than $B$ are rendered infeasible to the model~\eqref{prob:l2-FS-SVM2} and its relaxation~\eqref{prob:box-l2-FS-SVM1} below. 
\end{remark}

In the first formulation, we need a big-M parameter $M$ for the reformulation. It is known that the choice of $M$ has a profound influence on the computational performance if we solve the problem with off-the-shell global solvers. In a branch-and-bound framework, we solve LP relaxation at each node, and in many cases, the LP relaxation can be as weak as the SVM problem without feature selection, if too large $M$ is chosen, as shown later in Theorem~\ref{threshold_theorem}. In contrast, in the second formulation, the difficulty of the problem mainly goes to the complementarity constraint, and no big-M needed. The computational comparison between these two formulations will be presented in Sec. \ref{Sec5}. 

\section{Tractable Relaxations}\label{Sec3}
In this section, we propose several tractable relaxations, e.g., LP relaxations, and SDP relaxations, of the FS-SVM model, and a comparison is made among some of them. Moreover, a decomposition strategy is proposed by exploring the sparse patterns involved in  \red{both the} objective and the feasible region of the original problem, leading to a much smaller and scalable SDP-based relaxation. 

\subsection{Box relaxation}
In this subsection, we study the box relaxation (also known as a type of LP relaxation) of \red{the two} exact mixed-integer formulations presented in the last chapter and \red{highlight the crucial role that $M$ plays in the first formulation.}. 

After relaxing the binarity constraints to the box constraints in the problem \eqref{prob:l2-FS-SVM1}, we arrive at:

\begin{equation}
\label{prob:box-l2-FS-SVM1}
\begin{array}{rll}
    (\text{BoxMP}) \quad
    \min \limits_{\w,b,\bm{\xi},\vv } \quad &  \ds\frac 1 2 \|\w\|^2_2 + C \sum_{ i=1}^m \xi_{i} \\ 
    & ({\w,b,\bm{\xi} }) \in \Xi  \\[.8em] 
        & ({\w, \vv}) \in  \Omega_M \\[.8em] 
   & \ooo \leq \vv \leq \e\, .
    \end{array} 
\end{equation} 

The relaxation is a convex quadratic optimization problem solvable in polynomial time, therefore it has been widely used in many previous works, see references in~\cite{labbe2019mixed}. However, the tightness of the relaxation is highly dependent on the choice of $M$ and there is a possibility that it can be as weak as the original SVM model without feature selection~\eqref{prob:l2-SVM}. To be more precise, if $M$ is larger than a data-based threshold, the relaxation \eqref{prob:box-l2-FS-SVM1} is equivalent to the problem \eqref{prob:l2-SVM}, as presented in the following theorem. 

\begin{theorem}\label{threshold_theorem}
    If the parameter $M$ in the problem \eqref{prob:box-l2-FS-SVM1} satisfies $$M \geq \frac{\|\w^*\|_1}{B}\, ,$$ where $(\w^*,b^*, {\bm{\xi}}^*)$ is  an \eqref{prob:l2-SVM}-optimal solution, then \eqref{prob:box-l2-FS-SVM1} is equivalent to \eqref{prob:l2-SVM}.
\end{theorem}
\begin{proof}
    It is clear to see that \eqref{prob:l2-SVM} is a relaxation of \eqref{prob:box-l2-FS-SVM1}, as the former includes the constraint  $(\w , b,\bm\xi)\in \Xi$. Hence the optimal value of \eqref{prob:l2-SVM} is not larger than that of \eqref{prob:box-l2-FS-SVM1}. 
    To show the converse, we assume that $(\w^*,b^*,\bm\xi^*)$ is optimal to the problem~\eqref{prob:l2-SVM}. Let $v_j = \frac{|w^*_j|}{M}$ for all $j \irg{1}{n}$. Then the point $(\w^*,b^*,\bm\xi^*,\vv)$ is \eqref{prob:box-l2-FS-SVM1}-feasible. Indeed, due to the construction of $v_j$, the constraint $-M v_{j} \leq w_{j} \leq M v_{j}$ is satisfied automatically for all $j\irg{1}{n}$. Since $\e^\top \vv  = \frac{\|\w^*\|_1}{M}$, the constraint $\e^\top \vv \leq B$ is implied by the assumption we made. Since the objective functions in both problems coincide, the claim follows.
\end{proof}

On the other hand, $M$ can not be chosen smaller than: $$M \geq \|\overline{\w}^*\|_{\infty}\, ,$$ where $\overline{\w}^*$ is part of an optimal solution to~\eqref{prob:l2-FS-SVM}. Otherwise, the optimal solution will be cut off in the sense that the problem \eqref{prob:l2-FS-SVM1} is not a reformulation for an insufficient choice of $M$. Hence, the proper range of the parameter $M$ has to be:
\begin{equation}\label{bound_M}
    M\in\left[\|\overline{\w}^*\|_{\infty}, \frac{\|\w^*\|_1}{B}\right).
\end{equation}
However, this set can be empty depending on the problem data itself, which means no matter what $M$ we chose, the box relaxation~\eqref{prob:box-l2-FS-SVM1} ends up with one of the following situations: (i) no improvement compared to the original SVM problem~\eqref{prob:l2-SVM}; (ii) the optimal solution of the targeted problem \eqref{prob:l2-FS-SVM} is cut off, because $M$ is too small to render~\eqref{prob:box-l2-FS-SVM1} a reformulation of~\eqref{prob:l2-FS-SVM}. This is quite problematic in practice and it encourages us to explore the other formulation \eqref{prob:l2-FS-SVM2}.

\begin{remark}\label{rem3.2}
    Since estimating the range of $\overline{\w}^*$ is as difficult as solving the problem itself, we usually choose $M$ under a safe threshold, which is much larger than the tightest upper bound. In most practical cases, the condition $M \geq \frac{\|\w^*\|_1}{B}$ holds then as well, which means the LP relaxation~\eqref{prob:box-l2-FS-SVM1} does not provide any improvement compared to the plain problem~\eqref{prob:l2-SVM}.
\end{remark}

Due to the presence of complementarity constraints, the box relaxation of~\eqref{prob:l2-FS-SVM2} is equivalent 
to~\eqref{prob:l2-FS-SVM2} itself, as mentioned in Remark~\ref{remark2}. However, the combinatorial nature of a complementarity constraint is the main difficulty in both the unrelaxed and the relaxed version, so this continuous reformulation does not solve the problem immediately without any further processing.

\subsection{The Shor relaxation}

In this subsection,  we study the Shor relaxation, a type of SDP relaxation, of the two proposed mixed-integer reformulations.

Letting $\Wb = \w\w^\top $, $\Vb = \vv\vv^\top $ and dropping the rank one constraint, the Shor relaxation of \eqref{prob:l2-FS-SVM1} reads as:
\begin{equation}
\label{prob:sdp-l2-FS-SVM1}
\begin{array}{rll}
    (\text{SMP}) \quad
    \min \limits_{\w,\Wb,b,\bm{\xi } ,\vv,\Vb } \quad &  \ds\frac 1 2 \trace(\Wb) + C \sum_{ i=1}^m \xi_{i} \\ 
       \text{s.t.} \quad & ({\w,b,\bm{\xi } }) \in \Xi  \\[.8em]
       & ({\w,\vv}) \in \Omega_M  \\[.8em]
   & \diag \Vb =\vv  \\[.8em]
    & \begin{bmatrix}
        1 &\w^\top \\
        \w &\Wb \\
    \end{bmatrix} \in \SDP_+^{n+1},\ \begin{bmatrix}
        1 &\vv^\top \\
        \vv &\Vb \\
    \end{bmatrix} \in \SDP_+^{n+1} \, .
\end{array}
\end{equation}

\begin{theorem}
    The Shor relaxation~\eqref{prob:sdp-l2-FS-SVM1} is equivalent to the box relaxation \eqref{prob:box-l2-FS-SVM1}.
\end{theorem}

\begin{proof}
    For a \eqref{prob:sdp-l2-FS-SVM1}-feasible point $(\w,\Wb,b,\bm\xi ,\vv,\Vb)$, the projected point $(\w,b,\bm\xi,\vv)$ satisfies all of the constraints of \eqref{prob:box-l2-FS-SVM1}. Moreover, we have $\trace (\Wb) \ge \| \w\|_2^2$, hence the relaxation \eqref{prob:box-l2-FS-SVM1} is weaker than the relaxation \eqref{prob:sdp-l2-FS-SVM1}. 
    To show the other direction, we suppose that the point $(\w,b,\bm\xi,\vv)$ is feasible to the problem \eqref{prob:box-l2-FS-SVM1}, and let $\Wb = \w\w^\top $ and $\Vb = \vv\vv^\top + \Zb$, where $\Zb\in \SDP^n_+$ is 
    a diagonal matrix with $Z_{jj} = v_j - v_j^2 \geq 0$ for all $j\irg{1}{n}$. The constructed point $(\w, \Wb,b,\bm\xi,\vv,\Vb)$ then satisfies all constraints of~\eqref{prob:sdp-l2-FS-SVM1}. Indeed, the last positive-semidefiniteness constraint follows by the construction of $\Vb$ due to the Schur complement,
 which therefore completes the proof.
\end{proof}

Letting $\Wb = \w\w^\top $, $\Rb = \w\uu^\top $, $\Ub = \uu\uu^\top $ and dropping the rank-one constraint, we arrive at the Shor relaxation of \eqref{prob:l2-FS-SVM2}:
\begin{equation}
\label{prob:sdp-l2-FS-SVM2}
\begin{array}{rll}
    (\text{SCoP}) \quad
    \min \limits_{\substack{\w,\Wb,b,\bm{\xi}, \\ \Rb, \uu, \Ub} } \quad &  \ds\frac 1 2 \trace(\Wb) + C \sum_{ i=1}^m \xi_{i} \\ 
         \text{s.t.} \quad & ({\w,b,\bm{\xi } }) \in \Xi  \\[.8em]
    & \e^\top \uu = n - B  \\
    &\diag \Rb = \ooo, \quad 
     \diag \Ub =\uu
 \\[.8em]
    & \begin{bmatrix}
        1 &\w^\top &\uu^\top \\
        \w & \Wb & \Rb\\
        \uu & \Rb^\top & \Ub\\ 
    \end{bmatrix} \in \SDP_+^{2n+1}  \\
\end{array}\end{equation} 

\begin{theorem}
    If there exists part $\w$ of an optimal solution to~\eqref{prob:l2-SVM} with $\|\w\|_0 > B$, then the optimal value of~\eqref{prob:sdp-l2-FS-SVM2} is strictly larger than the optimal value of \eqref{prob:l2-SVM}.
\end{theorem}

\begin{proof}
   For a given \eqref{prob:l2-SVM}-feasible point $(\w,b,\bm\xi)$ with $\| \w\|_0 > B$, we restrict the corresponding variables of \eqref{prob:sdp-l2-FS-SVM2} to the same value and estimate the optimal value of the restricted version of \eqref{prob:sdp-l2-FS-SVM2}. The matrix 
$$\begin{bmatrix}
        1 &\w^\top &\uu^\top \\
        \w & \Wb & \Rb\\
        \uu &\Rb^\top & \Ub\\ 
    \end{bmatrix}
\in \SDP_+^{2n+1},\quad \mbox{with} \quad \diag \Rb = \ooo \mbox{ and }   \diag \Ub = \uu  \, ,$$ 
has a principal submatrix 
$$\begin{bmatrix}
        1 &w_j &u_j\\
        w_j &W_{jj} & 0\\
        u_j &0& u_j\\ 
    \end{bmatrix}\in \SDP_{+}^{3},\quad \mbox{for all }   j\irg{1}{n}\, ,$$ 
    which implies that 
    $$\left\{
    \begin{aligned}&0\leq u_j\leq 1\, , \quad W_{jj}\geq w_j^2, &&j\irg{1}{n}\\
    &(W_{jj} - w_j^2)(u_j - u_j^2)- w_j^2 u_j^2 \geq 0\, , &&j\irg{1}{n}\, .\\
    \end{aligned}\right. $$
    It follows that
    $$\left\{
    \begin{aligned}&0\leq u_j\leq 1\, , \quad &&\mbox{if} \quad j\in S_0 := \{j:w_j = 0\} \, ,\\
    & 0 \leq u_j\leq \frac{W_{jj}-w_j^2}{W_{jj}}<1 \quad&& \mbox{if} \quad j\in S_+ := [1\! : \!n] \setminus S_0\\
    \end{aligned}\right.  $$
    (observe that $w_j>0$ implies $W_{jj}>0$). Recalling  the assumption $\| \w\|_0 >B$, it follows that 
    $$\sum_{j\in S_0} u_j< n-B\, .$$
    Together with the constraint 
    $$\e^\top \uu = n-B\, ,$$
    we have $u_j > 0$ for some $j\in S_+$, which implies that $W_j > w_j^2$ for some $j\in S_+$. Hence, the optimal value of \eqref{prob:sdp-l2-FS-SVM2} is strictly larger than the optimal value of \eqref{prob:l2-SVM}. \end{proof}

There is no direct comparison between the relaxation \eqref{prob:sdp-l2-FS-SVM1} and the relaxation \eqref{prob:sdp-l2-FS-SVM2}. But given the case we have a valid estimate of $M$ at hand satisfying \eqref{bound_M}, the relaxation \eqref{prob:sdp-l2-FS-SVM2} can be tightened with a small computational cost by adding the big-M constraint via $ {\cal B}_M$, which reads as: 
\begin{equation}
\label{prob:sdp-l2-FS-SVM3}
\begin{array}{rll}
    (\text{SCoMP}) \quad
    \min \limits_{\substack{\w,\Wb,b,\bm{\xi }, \\ \Rb, \uu, \Ub}} \quad &  \ds\frac 1 2 \trace(\Wb) + C \sum_{ i=1}^m \xi_{i} \\ 
       \text{s.t.} \quad
    & {(\w,b,\bm{\xi } )} \in \Xi  \\[.8em]
        & {(\w,\uu)} \in  {\cal B}_M  \\[.8em]
    & \diag \Rb = \ooo\, , \quad \diag \Ub = \uu \, ,
    \\[.8em]
    & \begin{bmatrix}
        1 &\w^\top &\uu^\top \\
        \w & \Wb & \Rb\\
        \uu &\Rb^\top & \Ub\\ 
    \end{bmatrix} \in \SDP_+^{2n+1}\, .  \\
\end{array}
\end{equation}

\begin{theorem}
    The relaxation \eqref{prob:sdp-l2-FS-SVM3} is \red{at least as tight as} 
    the relaxations \eqref{prob:sdp-l2-FS-SVM1} and \eqref{prob:sdp-l2-FS-SVM2}.
\end{theorem}
\begin{proof}
It is clear that the relaxation \eqref{prob:sdp-l2-FS-SVM3} is tighter than the relaxation \eqref{prob:sdp-l2-FS-SVM2} due to the additional big-M constraint. To show the other claim, we construct a \eqref{prob:sdp-l2-FS-SVM1}-feasible point, given a \eqref{prob:sdp-l2-FS-SVM3}-feasible point. Suppose that the point $(\w,\Wb,b,\bm\xi, \Rb, \uu, \Ub)$ is \eqref{prob:sdp-l2-FS-SVM3}-feasible, let $\vv = \e - \uu$ and  $\Vb = \vv\vv^\top + \Zb$, where $\Zb\in \SDP^n_+$ is a diagonal matrix with $Z_{jj} = v_j - v_j^2 \geq 0$ for all $j\irg{1}{n}$. Then $(\w,\Wb,b,\bm\xi, \vv, \Vb)$ is \eqref{prob:sdp-l2-FS-SVM1}-feasible with the same objective value. 
\end{proof}

\subsection{Decomposed SDP-based relaxations}

In this subsection, we propose a procedure for obtaining equivalent and much smaller SDP-based relaxations by exploiting the sparsity patterns of the Shor relaxations.

 We notice that, for example: in problem~\eqref{prob:sdp-l2-FS-SVM2}, some sparsity patterns appear both in the objective and the constraints. e.g., only diagonal elements of matrices $\Wb$, $\Rb$, and $\Ub$ are involved. Based on this observation, we write down the following decomposed version of~\eqref{prob:sdp-l2-FS-SVM2}:
 
\begin{equation}
\label{prob:dsdp-l2-FS-SVM2}
\begin{array}{rll}
    (\text{DSCoP}) \quad
     \min \limits_{{\w,\begin{footnotesize}\diag\end{footnotesize}\Wb,b,\bm\xi ,\uu} } \quad &  \ds\frac 1 2 \trace(\Wb) + C \sum_{ i=1}^m \xi_{i} \\ 
             \text{s.t.} \quad
    & {(\w,b,\bm{\xi } )} \in \Xi  \\[.8em]
       & \e^\top \uu = n - B \\[.8em]
    & \begin{bmatrix}
        1 &w_j &u_j\\
        w_j &W_{jj} &0\\
        u_j &0 &u_j\\ 
    \end{bmatrix} \in \SDP_+^3 \, , \, & j\irg{1}{n}   \\
\end{array} \end{equation}

\begin{theorem}\label{decomposed_theorem}
    The relaxation \eqref{prob:sdp-l2-FS-SVM2} is equivalent to the relaxation \eqref{prob:dsdp-l2-FS-SVM2}.
\end{theorem}

\begin{proof}
 Suppose that $(\w,\Wb,b,\bm\xi, \Rb, \uu, \Ub)$ is \eqref{prob:sdp-l2-FS-SVM2}-feasible. It is clear that $(\w,\diag\Wb,b,\bm\xi, \uu)$ is \eqref{prob:dsdp-l2-FS-SVM2}-feasible. In the opposite direction, we assume that the point $(\w,\diag \Wb,b,\bm\xi, \uu)$ is \eqref{prob:dsdp-l2-FS-SVM2}-feasible. It is clear that all linear constraints in \eqref{prob:sdp-l2-FS-SVM2} are satisfied if we let $\diag \Rb  = \ooo$ and $\diag \Ub  = \uu$. The only thing that needs to be checked is if the following partial positive-semidefinite matrix is positive-semidefinite completable:
$$\begin{bmatrix}
    1 & w_1 & \dots & w_n & u_1 & \dots & u_n\\
    w_1 & W_{11} & * & *& 0 & * & *\\
    \vdots & * & \ddots & *&* & \ddots & *\\
    w_n & * & * & W_{nn} & * & * & 0\\
    u_1 & 0 & * & *& u_1 & * & *\\
    \vdots & * & \ddots & *&* & \ddots & *\\
    u_n & * & * & 0 & * & * & u_n\\
\end{bmatrix}\, .$$
By $\mathcal{G}$ we denote the specification graph of the above partial positive-semidefinite matrix. Since the vertex $w_i$ is not connected with $w_j$ and $u_j$ for all $i\neq j$, and the vertex $u_i$ is not connected with $u_j$ and $w_j$ for all $i\neq j$, the only cycle we can find in the graph is the triangle $1\to w_{j} \to u_{j} \to 1$ for all $[1\!:\!n]$. Therefore $\mathcal{G}$ is chordal. Due to \cite[Theorem~1.39]{berman_completely_2003}, if the specification graph of a partial positive-semidefinite matrix is chordal, it is always completable to a positive-semidefinite matrix. 
\end{proof}

Following a similar fashion, we propose the decomposed version of the relaxation \eqref{prob:sdp-l2-FS-SVM3}:
\begin{equation}
\label{prob:dsdp-l2-FS-SVM3}
\begin{array}{rll}
    (\text{DSCoMP}) \quad
     \min \limits_{{\w,\begin{footnotesize}\diag\end{footnotesize}\Wb,b,\bm{\xi }  ,\uu} } \quad &  \ds\frac 1 2 \trace(\Wb) + C \sum_{ i=1}^m \xi_{i} \\ 
       \text{s.t.} \quad
   & {(\w,b,\bm{\xi }  )} \in \Xi  \\[.8em]
            & {(\w,\uu)} \in  {\cal B}_M  \\[.8em]
  & \begin{bmatrix}
        1 &w_j &u_j\\
        w_j &W_{jj} &0\\
        u_j &0 &u_j\\ 
    \end{bmatrix} \in \SDP_+^3\, , \,  & j\irg{1}{n} \, .\\
\end{array}
\end{equation}

\begin{corollary}
    The relaxation \eqref{prob:dsdp-l2-FS-SVM3} is equivalent to the relaxation \eqref{prob:sdp-l2-FS-SVM3}.
\end{corollary}

\begin{proof}
    See the proof of Theorem \ref{decomposed_theorem}.
\end{proof}

For the formulation \eqref{prob:l2-FS-SVM1}, the decomposed version of its Shor relaxation~\red{\eqref{prob:sdp-l2-FS-SVM1}}

reads as:

\begin{equation}
\label{prob:dsdp-l2-FS-SVM1}
\begin{array}{rll}
    (\text{DSMP}) \quad
     \min \limits_{{\w,\begin{footnotesize}\diag\end{footnotesize}\Wb,b,\bm{\xi} ,\vv} } \quad &  \ds\frac 1 2 \trace(\Wb) + C \sum_{i=1}^m \xi_{i} \\ 
       \text{s.t.} \quad
   & ({\w,b,\bm{\xi} }) \in \Xi  \\[.8em]
       & ({\w,\vv}) \in \Omega_M  \\[.8em]
    & \begin{bmatrix}
        1 &w_j\\
        w_j &W_{jj}\\ 
    \end{bmatrix} \in \SDP_+^2\, ,\,  
    \begin{bmatrix}
        1 &v_j\\
        v_j &v_{j}\\ 
    \end{bmatrix} \in \SDP_+^2\, , \, & j\irg{1}{n} \, .  \\
 \end{array}
\end{equation}

\begin{corollary}
    The relaxation \eqref{prob:dsdp-l2-FS-SVM1} is equivalent to the relaxation \eqref{prob:sdp-l2-FS-SVM1}.
\end{corollary}

\begin{proof}
    See the proof of Theorem \ref{decomposed_theorem}.
\end{proof}

\section{The algorithms}\label{Sec4}

In this section, we introduce both a heuristic and an exact algorithm for solving the feature selection problem we are interested in. In particular, we first present two heuristic strategies to find good upper bound values for our problem and a heuristic strategy to estimate a tight value for the big-M parameter $M$ for the problem. All strategies are based on the resolution of the decomposed SDP-based relaxations presented in the previous section. These strategies will then be embedded in a heuristic algorithm which can ultimately improve the heuristic upper bound. Finally, we present an exact algorithm that implements one of the two upper bound procedures as a subroutine and consists of the resolution of a sequence of Mixed-Integer Second-Order Cone Optimization problems (MISOCPs). 

\subsection{Two upper bounding strategies}

In this subsection, we propose two strategies for obtaining an upper bound, namely a feasible objective value of a feasible solution $(\widehat{\w},\widehat{b},\widehat{\bm\xi})$ to~\eqref{prob:l2-FS-SVM}, based on the optimal solution of relaxations. 

\subsubsection{Local Search}
The first strategy, denoted as Local Search, consists of three main steps and employs a user-specified excess parameter $k$ which will temporarily allow to work on just $k+B$ features instead of all. We first solve the relaxation by off-the-shelf software such as $\texttt{Mosek}$, an interior point solver. With the optimal solution $\uu^*$ in hand, we want to search for a feasible point nearby. We \red{then} sort features \red{by increasing values of} $\uu^*$. Namely, if $u^*_j$ is close to 1, feature $j$ is unlikely to be selected because of the complementarity constraint $w_j u_j = 0$, while if $u^*_j$ is closer to 0, then the feature $j$ is more likely to be selected. Based on it, we select the first $k+B$ features with small corresponding $u$, denoted by a set $K\subseteq [1\!:\!n]$ with $|K|= k+B$. In the end, we solve the restricted mixed-integer problem:
\begin{equation}
\label{prob:l2-FS-SVM2(K)}
\begin{array}{rll}
    (\text{CoP($K$)}) \quad
    \min \limits_{{\w,b,\bm\xi ,\uu} } \quad &  \ds\frac 1 2 \| \w\|^2_2 + C \sum_{ i=1}^m \xi_{i} \\ 
       \text{s.t.} \quad
    &({\w,b,\bm{\xi} }) \in \Xi  \\[.8em]
    & \sum\limits_{j\in K} u_j = \red{|K|} - B \\[.8em]
    &u_{j} w_{j} = 0 \, , &  j\in K \, , \\[.8em]
    & u_{j} \in\{0,1\} \, , & j\in K \, .
\end{array} \end{equation}
which is much smaller than the original problem \eqref{prob:l2-FS-SVM2}, as it only involves $|K|$ binary variables $u_j$. The problem is solved by the global solver $\texttt{Gurobi}$. This first upper bounding strategy is described in Algorithm~\ref{ub_algorithm_1}.

\begin{algorithm}\caption{ \texttt{Local Search} }\label{ub_algorithm_1}
\begin{algorithmic}[1]
    
\State {\bf Input}: ($\x^{i}$,$y^{i}$)$\in \mathbb{R}^{n}\times \{1,-1\}$, $B$, $C$, parameter $k \irg{0}{n-B}$ 
\State Solve the relaxation \eqref{prob:dsdp-l2-FS-SVM2} and return the optimal solution $\uu^*$
\State  Generate  a subset $K\subseteq [1\!:\!n]$ with $|K| = k+B$ such that  $u^*_i\leq u^*_j$ for all $i\in K$ and for all $j\in [1\!:\!n]\setminus K$
\State  Solve the mixed-integer problem \eqref{prob:l2-FS-SVM2(K)}
\State  {\bf Output}: The set of selected features and an upper bound UB.
 \end{algorithmic}
\end{algorithm}

\subsubsection{Kernel Search}

Slightly enlarging the search set at each iteration, we propose the following search strategy for a feasible solution of \eqref{prob:l2-FS-SVM}.
Our Kernel Search strategy consists in solving a sequence of small MIQPs considering an initial ranking on the features. This strategy takes inspiration from a heuristic proposed by \cite{angelelli2010kernel} for the general multi-dimensional knapsack problem which has been applied to different kinds of problems such as location problems \cite{guastaroba2012kernel} and portfolio optimization \cite{angelelli2012kernel}. Recently, \cite{labbe2019mixed} also applied it to solve their embedded feature selection for $\ell_1$-SVMs. In this strategy, based on some ranking, a restricted set of promising features is kept and updated throughout each iteration.
Usually, for ranking variables in general MILPs, LP theory is exploited using both relaxed solution values and reduced costs, a strategy followed by~\cite{labbe2019mixed} for ranking features. 

By contrast, in our case we decided to use as a ranking criterion the coordinates of an optimal solution  $\uu^*$ to the decomposed SDP relaxation~\eqref{prob:dsdp-l2-FS-SVM2}. Features are sorted  
with respect to this vector $\uu^*$: the smaller the value of $u_j^*$, the higher the relevance of feature $j$.
Since in practice~\eqref{prob:dsdp-l2-FS-SVM2} leads to much tighter lower bound solutions in contrast to the plain relaxation~\eqref{prob:box-l2-FS-SVM1} similar to the model used by~\cite{labbe2019mixed}, there is a justified hope that the ranking provided by $\uu^*$ is closer to the relevance of features.  

Let $K$ be the ordered set of features,  $\rho$ a user-defined parameter and $N:=\lceil \frac{n-\rho}{\rho}\rceil$. Then $K$ is divided into $N$ subsets denoted as $K_i$ for $i \irg 1N $.
In particular, each subset $K_i$, $i\irg 1{N-1} $, will be composed of $\rho$ features, and $K_N$ will contain the remaining features. 

An initial value of the UB is set to $\infty$. Also define a feature subset $\overline K \subseteq [1\!:\!n]$ as the \textit{kernel set}, containing features that are kept at the next iteration. To initialize, set $\overline K=\emptyset$. At each iteration $k$, the heuristic considers set  $\mathcal{K}_k = \overline K \cup K_k$ i.e. the union of the features in the kernel set and the features in the set $K_{k}$. To update the UB, at each iteration we solve problem \text{CoP($\mathcal{K}_k$)} plus the following two constraints:

\begin{equation}
\ds\frac 1 2 \|\w\|^2_2 + C \sum_{ i=1}^m \xi_{i} \leq \mbox{UB}, \label{cons: UB cons}
\end{equation}
\begin{equation}
    \ds\sum_{j \in K_{k}} u_{j} \leq |K_{k}| - 1.
    \label{cons: K plus cons}
\end{equation}

Constraint \eqref{cons: UB cons} restricts the objective function to take a value smaller than or equal to the current upper bound, while constraint \eqref{cons: K plus cons} makes sure that at least one feature from the new set $K_{k}$ is selected. Note that, due to the addition of these constraints, the optimization problems may potentially be infeasible. If the optimization problems are feasible, the new features selected are then added to the kernel set $\overline K$ used in the next iteration since adding these features obtains an identical or better upper bound. Conversely, the set
of features of $\overline K$ that have not been chosen in the optimal solution in the previous iterations is removed from the kernel set $\overline K$. The
removal of some of the features from the kernel set is decisive in that it does not excessively increase the number of binary variables considered in each iteration. We decided to remove the features that were not selected in the
previous two iterations. The set of added features is denoted as $K^+$ and the set of removed features as $K^-$. The resulting
kernel set for the next iteration is $\overline{K} \cup K^{+} {\setminus  K^{-} }$.
Conversely, if the problem is infeasible, the kernel set is not modified and the procedure skips to the next iteration. The Kernel Search strategy is described in Algorithm \ref{ub_algorithm_2}.

\begin{algorithm}\caption{ \texttt{Kernel Search} }\label{ub_algorithm_2}
 \begin{algorithmic}[1]
\State {\bf Input}: ($\x^{i}$,$y^{i}$)$\in \mathbb{R}^{n}\times \{1,-1\}$, $B$, $C$, parameter $\rho \irg 1n$ 
\State Solve the relaxation \eqref{prob:dsdp-l2-FS-SVM2} and return the optimal solution $\uu^*$
\State Sort features according to the numerical size of corresponding $\uu^*$ increasingly, denoted by a vector $K\in \mathbb{Z}_+^n$ 
\State Separate $K$ into $N:=\lfloor\frac{n}{\rho}\rfloor$ groups evenly such that $K = [K_1, \cdots, K_N]$ and set $\overline{K} = \emptyset$,  $ k = 1$, UB $= \infty$
\While{$k \leq N$}
    \State $\mathcal{K}_{k} \gets K_{k} \cup \overline{K}$
    \State Solve the problem CoP($\mathcal{K}_{k}$) + \eqref{cons: UB cons} + \eqref{cons: K plus cons} \If{problem is feasible}
    \State Return  solution  $\uu^*$ and optimal value UB$^*$
    \State UB $\gets$ UB$^*$
    \State Build $K^+:= \{j \in K_{k}: j \text{ selected at current iteration}\}$
    \State Build $K^-:= \{j \in \overline K: j \text{ not selected in the solution of the last two iterations}\}$
    \State Update $\overline{K} \gets  {\overline K\cup K^{+} }{\setminus  K^{-} }$
\EndIf
 \State $k \gets k + 1$
\EndWhile
\State  {\bf Output}: The set of selected features and an upper bound UB.
\end{algorithmic}
\end{algorithm}

\subsection{A strategy to tighten the big-M parameter}

Given an upper bound UB of problem~\eqref{prob:l2-FS-SVM}, one can estimate the upper bound of $w_j$ by solving the following problem:

\begin{equation}
\label{subproblem_bound_M}
\begin{array}{rll}
    \quad
    \overline{M}_j := \max \limits_{{\w,\begin{footnotesize}\diag\end{footnotesize}\Wb,b,\bm\xi ,\uu} } \quad & w_j \\ 
       \text{s.t.} \quad
     &({\w,b,\bm{\xi} }) \in \Xi  \\[.8em]
    & \ds\frac 1 2 \|\w\|^2_2 + C \sum_{ i=1}^m \xi_{i} \leq \mbox{UB}  \\[.8em]
      & \e^\top \uu = n - B \\[.8em]
    & \begin{bmatrix}
        1 &w_j &u_j\\
        w_j &W_{jj} &0\\
        u_j &0 &u_j\\ 
    \end{bmatrix} \in \SDP_+^3\, , \, & j\irg{1}{n} \, .
\end{array}\end{equation}

Similarly, the lower bound of $w_j$ can be obtained by solving the same optimization problem with the max replaced by min and we denote the optimal value by ${\underline{M}}_j$. The proper choice of $M$ is $\max\limits_{j\irg{1}{n}}\{\max\{|\overline{M}_j|,|\underline{M}_j|\}\}$.
\begin{proposition}\label{prop_M}
    Let $[\widehat{\w},\widehat{b},\widehat{\bm\xi}]$ be a feasible solution of \eqref{prob:l2-FS-SVM} and $\mbox{UB} :=  \ds\frac 1 2 \|\widehat{\w}\|^2_2 + C \sum_{ i=1}^m \widehat{\xi}_{i}$ be the corresponding objective value. Then 
    $$M = \max_{j\in [1:n]}\{\max\{|\overline{M}_j|,|\underline{M}_j|\}\}$$ 
is a valid bound, where $\overline{M}_j$ and $\underline{M}_j$ are obtained 
by~\eqref{subproblem_bound_M}.
\end{proposition}
\begin{proof} obvious. \end{proof}

\begin{remark}
    The proposed bound $M$ is tighter than most of the currently used big-Ms.  
    Nevertheless, empirically, the tightening of any such bound is of little relevance unfortunately, because they still fall outside of the range specified in~\eqref{bound_M}.
\end{remark}

\subsection{The heuristic algorithm}
In this subsection, we propose a procedure for generating a good feasible solution of \eqref{prob:l2-FS-SVM} taking one of the two upper bound strategies and the big-M estimation strategy into account. At first, a feasible solution is computed with the chosen upper bound strategy. Then, according to Proposition \ref{prop_M}, we solve a series of 2$n$ many SDP problems to find, for each feature $j$, values $\overline{M}_j$,$\underline{M}_j$ and then obtain the big-M parameter $M$.  If this value satisfies \eqref{bound_M}, we know that the SDP formulation \eqref{prob:dsdp-l2-FS-SVM3} is better than \eqref{prob:dsdp-l2-FS-SVM2}, thus we compute a new relaxed solution $\uu^*$ of \eqref{prob:dsdp-l2-FS-SVM3}. We thus use vector $\uu^*$ as a ranking criterion for the chosen upper bound strategy, we run the strategy again and we update the heuristic solution if a better one is found.
\begin{algorithm}\caption{ \texttt{Heuristic procedure}}\label{ub_algorithm_3}
\begin{algorithmic}[1]
\State {\bf Input}: ($\x^{i}$,$y^{i}$)$\in \mathbb{R}^{n}\times \{1,-1\}$, $B$, $C$, the parameter $k \irg 0{n-B}$ or the parameter $\rho \irg 1n$ 
\State Run Algorithm~\ref{ub_algorithm_1} or Algorithm~\ref{ub_algorithm_2} to obtain a feasible solution $(\widehat{\w},\widehat{b},\widehat{\bm{\xi}})$ to~\eqref{prob:l2-FS-SVM} with upper bound value $\widehat{\text{UB}}$
\State Set $(\w^*,b^*,\bm\xi^*) \gets (\widehat{\w},\widehat{b},\widehat{\bm\xi})$
\State Obtain $M$ according to Proposition~\ref{prop_M}
\State Solve problem \eqref{prob:l2-SVM} and obtain $\w^*$
    \If{$M < \frac{\|\w^*\|_1}{B}$}
        \State Solve problem \eqref{prob:dsdp-l2-FS-SVM3} and obtain $\uu^*$
        \State Run  Algorithm~\ref{ub_algorithm_1} or Algorithm~\ref{ub_algorithm_2} starting with $\uu^*$, get $(\overline {\w},\overline{b},\overline{\bm{\xi}})$ feasible to~\eqref{prob:l2-FS-SVM} and upper bound $\overline{\text{UB}}$
        \If{$ \overline{\text{UB}} \leq \widehat{\text{UB}}$}
        \State $(\w^*,b^*,\bm\xi^*) \gets (\overline{\w},\overline{b},\overline{\bm\xi})$
                \EndIf
        \EndIf 
\State {\bf Output} $(\w^*,b^*,\bm\xi^*)$ solution of  \eqref{prob:l2-FS-SVM}
\end{algorithmic}
\end{algorithm}

\subsection{The exact algorithm}
In this subsection, we propose an exact procedure for \eqref{prob:l2-FS-SVM} by solving a sequence of MISOCPs. With a slight abuse of notation, let us  define, given a subset of features $K\subseteq [1\!:\!n]$, the set
\begin{equation*}
{\cal B}_M(K)=\{(\w, \uu )\in\red{\R^n\times \R^n_+}: -M(1-u_j)\leq w_j \leq M(1-u_j),\, j\in K\,, \ \e^\top \uu = n - B \}\, .    
\end{equation*} The subproblem considered in the exact procedure is

\begin{equation}
 \label{subproblem_exact_procedure}
\begin{array}{rll}
    (\text{SR-DSCoP($K$)}) \quad
    \min \limits_{{\w,\begin{footnotesize}\diag\end{footnotesize}\Wb,b,\bm\xi ,\uu} } \quad &  \ds\frac 1 2 \trace(\Wb) + C \sum_{ i=1}^m \xi_{i} \\[.8em] 
       \text{s.t.} \quad
    & (\w,b,\bm\xi)\in {\Xi}\\[.8em]
    & (\w, \uu ) \in {\cal B}_M(K)\\[.8em]
    & \begin{bmatrix}
        1 &w_j &u_j\\
        w_j &W_{jj} &0\\
        u_j &0 &u_j\\ 
    \end{bmatrix} \in \SDP_+^3 \, , & j\irg{1}{n}\setminus  K \,  , \\[.8em]
    & u_j\in \{0,1\} \, , \, & j \in K\, ,  \\[.8em]
     & \begin{bmatrix}
        1 &w_j\\
        w_j &W_{jj}\\
    \end{bmatrix} \in \SDP_+^2 \,  , \, & j \in K \, , 
\end{array}
\end{equation}
where SR is the abbreviation for semi-relaxation. In the problem \eqref{subproblem_exact_procedure}, we essentially relax the complementarity constraints associated with features in the set $[1\!:\!n]\setminus K$ via decomposed SDP relaxation, and express the other by the big-M mixed-integer formulation. We note that the problem~\eqref{subproblem_exact_procedure} is a mixed-integer positive semi-definite optimization problem (MISDP)  and there is no off-the-shelf solver dealing with it. In the following part, we show a way to express the problem in the form of MISOCP, which can be handled by the off-the-shelf solvers like~$\texttt{copt}$.

Before specifying the MISOCP reformulation, we present a useful result: some special slices $\mathcal{F}$  of the positive-semidefinite cone  are second-order cone representable, due to Lemma~\ref{PSD_SOC} below
\red{which seems well known in the literature, see, e.g.~\cite{agler1988positive,cheramin2022computationally}. We include its short, and elementary, proof for the sake of completeness.}

\begin{lemma}\label{PSD_SOC}
 \begin{equation*}
         \mathcal{F}:=\left\{(w,u,W)\in \R^3: \begin{bmatrix}
             1 & w & u\\
             w & W & 0\\
             u & 0 & u\\
         \end{bmatrix}\in \SDP_+^3 \right\} = 
         \left\{(w,u,W)\in \red{\R\times \R^2_+}:
         \red{\begin{aligned}
             &\sqrt{(W-s)^2 + 4 w^2} \leq W+s \\
                          &s := 1-u \ge 0\, .
         \end{aligned} } 
         \right\}
 \end{equation*}   
\end{lemma}

\begin{proof}
   \red{Consider all principal minors of the matrix (with a chordal sparsity pattern, see~\cite{agler1988positive}):} 
    $$\begin{bmatrix}
             1 & w & u\\
             w & W & 0\\
             u & 0 & u\\
         \end{bmatrix}\in \SDP_+^3 \quad \Longleftrightarrow\quad 
   \red{(u,W)\in \R^2_+ \quad\mbox{and} \quad u\le 1 \quad\mbox{and} \quad  (1-u)W\ge w^2\, . }   $$

     \red{Putting $s=1-u$,} it further follows that \red{for $(s,W)\in \R^2_+ $, we have $sW\ge w^2 \Longleftrightarrow 
     \sqrt{(W-s)^2 + 4 w^2} \leq W+s$.
          This way}
         we complete the proof.
\end{proof}

Applying Lemma~\ref{PSD_SOC}, the problem \eqref{subproblem_exact_procedure} can be rewritten of the form
\begin{equation}
\label{subproblem_exact_procedure_modified}
\begin{array}{rll}
    (\text{SR-DLMP($K$)}) \quad
    \min \limits_{\substack{\w,\begin{footnotesize}\diag\end{footnotesize}\Wb,b,\bm\xi ,\uu,\sss} } \quad &  \ds\frac 1 2 \trace(\Wb) + C \sum_{ i=1}^m \xi_{i} \\[0.8em] 
       \text{s.t.} \quad
    & (\w,b,\bm\xi)\in {\Xi}\\[0.8em]
    & (\w, \uu ) \in {\cal B}_M(K)\\[0.8em]
    & \ddd_j :=\begin{bmatrix}
        W_{jj} - s_j \\
        2 w_j
    \end{bmatrix}\in \R^2\, ,  \,  & j\irg{1}{n}\setminus  K  \, ,\\[0.8em]
    & \|\ddd_j\|_2\leq s_j + W_{jj}\, , \,  &j\irg{1}{n}\setminus  K \, ,  \\[0.8em]
    & s_j + \red{u_j} = 1 \, , \, & j\irg{1}{n}\setminus  K\, ,  \\[0.8em]
    & u_j \in \{0,1\} \, , \, & j \in K\, , \\[0.8em]
    & w_j^2 \leq W_{jj}\, , \,  & j \in K\, ,  \\[0.8em]
\end{array}\end{equation}   
where \red{$n-|K|$  more variables $s_j$, $j\irg{1}{n}\setminus  K$,
are introduced.}

Now we can introduce the exact algorithm reported as Algorithm \ref{algo: sfs-e}. In this exact procedure, we solve a sequence of semi-relaxed problems SR-DLMP($K$) associated with a subset $K \subseteq [1\!:\!n]$ of features. In such problems, only variables $u_j$ with $j \in K$ will be considered as binary and the remaining ones are relaxed. To obtain initial bounds on the objective value, the exact procedure exploits the upper bound strategies detailed above. The main step of the exact procedure consists of solving the sequence of semi-relaxed problems~\eqref{subproblem_exact_procedure_modified} to improve the lower bound of the objective value. 

To start with, we must select a subset of features $K\subseteq [1\!:\! n]$  whose associated $\uu$ variables will be considered as binary variables in the first semi-relaxed problem. The upper-bound strategy run at the beginning provides a subset of features that allows us to obtain a good bound on the optimal objective value. Therefore, the exact procedure will consider the set provided by the heuristic as the initial $K$ and it will obtain an initial LB solving SR-DLMP($K$). Then the set $K$ is updated by adding and removing some of the features, in order to improve the lower bound of the objective value. To this end, two sets (denoted by $K^+$ and $K^-$) are built in each iteration. Set $K^+$ consists of some of the features in $  [1\! : \! n]\setminus K$ whose associated $\uu$ variables will be considered as binary in the next iteration, i.e. features of $K^+$ will be added to $K$. Similarly, $K^-$ consists of features in $K$ that will not be considered as binary in the next iteration. 

\red{While we cannot rule  out that a once selected set $K$ is visited again during the iterations, the imposed time limit ensures finiteness of this procedure.}

In addition and if possible, we will update the UB in the main step running the upper bound strategies and using vector $\widetilde \uu$ as the initial ranking. There could be many different ways of defining set $K^+$ and updating $K$. In the strategy we adopted, we set $K^+$  as the set of $s$ features $j$ such that the relaxed solution $ \red{\widetilde u_j}$ of SR-DLMP($K$) was "less binary". This is why we first create a ranking on the features based on terms $\widetilde u_j - \widetilde u^2_j$ in descending order, and then select the first $s$, which is a user-defined parameter.

\begin{algorithm}
\caption{\texttt{Exact procedure}}
\label{algo: sfs-e}
\begin{algorithmic}[1]
\State \textbf{Input}: $(\x^{i},y^{i})\in \mathbb{R}^{n}\times \{1,-1\}$, $B$, $C$,  $s \irg{1}{n}$, LB = $-\infty$, UB = $\infty$, $K = \emptyset$
\State Run Algorithm~\ref{ub_algorithm_1} or Algorithm~\ref{ub_algorithm_2}, obtain set $K$ as the set of the selected features and update the UB value
\While{$\frac{UB-LB}{UB} \cdot 100 \geq 0.01$}
    \State Solve problem SR-DLMP($K$)
    \If{a time limit of $1800$ seconds is reached}
    \State {\bf STOP}
    \EndIf
    \State Obtain optimal value $\widehat{\text{LB}}$ and optimal solution $\widetilde \uu$
    \State Run Algorithm~\ref{ub_algorithm_1} or Algorithm~\ref{ub_algorithm_2} with the first step skipped vector $\widetilde \uu$ used for initial ranking
    \State Obtain the solution $(\widehat \w,\widehat b,\widehat{\bm{\xi}} , \widehat \uu)$ and upper bound value $\widehat{\text{UB}}$
    \State Sort features using terms $\widetilde u_j - \widetilde u^2_j$ in descending order 
    \State Build set $K^+$ with the first $s$ features
    \State Build $K^-:= \{j \in K: j \text{ not selected in the solution of the last two iterations}\}$
    \State $K \gets K \cup K^+ \setminus K^-$
\If{$\widehat{\text{UB}} < \text{UB}$}
    \State UB $\gets$ $\widehat{\text{UB}}$
    \State $(\w^*,b^*,{\bm{\xi}}^* ,{\uu}^*) \gets (\widehat \w,\widehat b,\widehat {\bm{\xi}} , \widehat \uu)$
\EndIf
\If{$\widehat{\text{LB}}>$ LB}
    \State LB $\gets$ $\widehat{\text{LB}}$
\EndIf
\EndWhile
\State {\bf Output}:   $(\w^*,b^*,{\bm{\xi}}^*)$ solution of  \eqref{prob:l2-FS-SVM}
\end{algorithmic}

\end{algorithm}

\section{Numerical experiments}\label{Sec5}

\begin{table}[ht]
\centering
\begin{minipage}[b]{0.45\linewidth} 
\centering
\renewcommand\arraystretch{1.2}
\begin{tabular}{lcccc}
\hline
Dataset                  &  & $m$   & $n$ & Class (\%) \\ \hline \hline
Breast Cancer D. &  & 569 & 30  & 63/37      \\
Breast Cancer P.  &  & 198 & 33   & 65/35      \\
Breast Cancer W.  &  & 683 & 9   & 65/35      \\
Cleveland  &  & 297 & 13  & 54/46      \\
Diabetes  &  &  768 & 8  & 65/35     \\
German &  &  1000 & 24  & 70/30      \\
Ionosphere        &  & 351 & 33  & 36/64      \\
Parkinsons         &  & 195 & 22  & 25/75      \\
Wholesale                &  & 440 & 7   & 68/32   \\ \hline  
\end{tabular}
\caption{First \red{small} datasets}
\label{table:small_datasets}
\end{minipage}
\hfill
\begin{minipage}[b]{0.45\linewidth} 
\centering
\renewcommand\arraystretch{1.2}
\begin{tabular}{lcccc}
\hline
Dataset                  &  & $m$   & $n$ & Class (\%) \\ \hline \hline
Arrhythmia &  & 452 &  274  &   54/46    \\
Colorectal &  & 62 & 2000   &  65/35   \\
DLBCL  &  & 77 & 7129  &   25/75    \\
Lymphoma &  &  96 & 4026  &  35/65  \\
Madelon               &  & 2000 & 500  & 50/50     \\
Mfeat               &  &2000  & 649   & 10/90    \\ \hline
\end{tabular}
\caption{Second \red{large} datasets}
\label{table:big_datasets}
\end{minipage}
\end{table}

The various numerical experiments were performed on an Apple M1 CPU 16 GB RAM computer. All models and algorithms were implemented in \texttt{Python}.  
Results were performed using \texttt{Gurobi 11.0} for the resolution of the MIQPs, \texttt{Mosek 10.1} for the resolution of the SDP relaxations and \texttt{Copt} for the resolution of the MISOCP models solved in the Exact algorithm. 

The experiments were carried out on fifteen different datasets. Nine of them can be found in the UCI repository \cite{UCI2019}. See Table \ref{table:small_datasets}), where $m$ is the number of data points, $n$ is the number of features and the last column shows the percentage of samples in each class. As can be observed, they contain a small number of features. The other six datasets used in the experiments have a larger number of features, as shown in Table \ref{table:big_datasets}. The Arrythmia, Madelon, and MFeat datasets are available in the UCI repository as well. The remaining three datasets, Colorectal, DLBCL, and Lymphoma, are microarray datasets with thousands of features but smaller sample sizes. Further descriptions of these last datasets can be found in~\cite{alon1999broad,guyon2002gene,shipp2002diffuse}.

\subsection{Computational analysis on the MIQP models}

\red{To start with, we used \texttt{Gurobi} {\cite{gurobi}} to solve} the two proposed MIQP models BigMP and CoP. 
\red{The reason behind these first comparisons is that \texttt{Gurobi} is used in the two upper bounding strategies we proposed. The two strategies, the Local Search and the Kernel Search, require the resolution of the FS-SVM model on instances created on a smaller subset of features, and this can be either done by using the BigMP or the CoP formulation.}

\red{
From these initial tests, we understand that CoP is the model which can be solved faster with \texttt{Gurobi} when dealing with datasets with a small number of features. Thus  we will use the resolution of model CoP as the subroutine used by the heuristic algorithms to find, and update the values of upper bounds.}

\red{We designed two sets of experiments: the first using small datasets, described in} Table~\ref{table:small_datasets} \red{and the second using larger ones, described in} Table~\ref{table:big_datasets}.

\begin{figure}[!h]
        \centering
        \includegraphics[width=18cm]{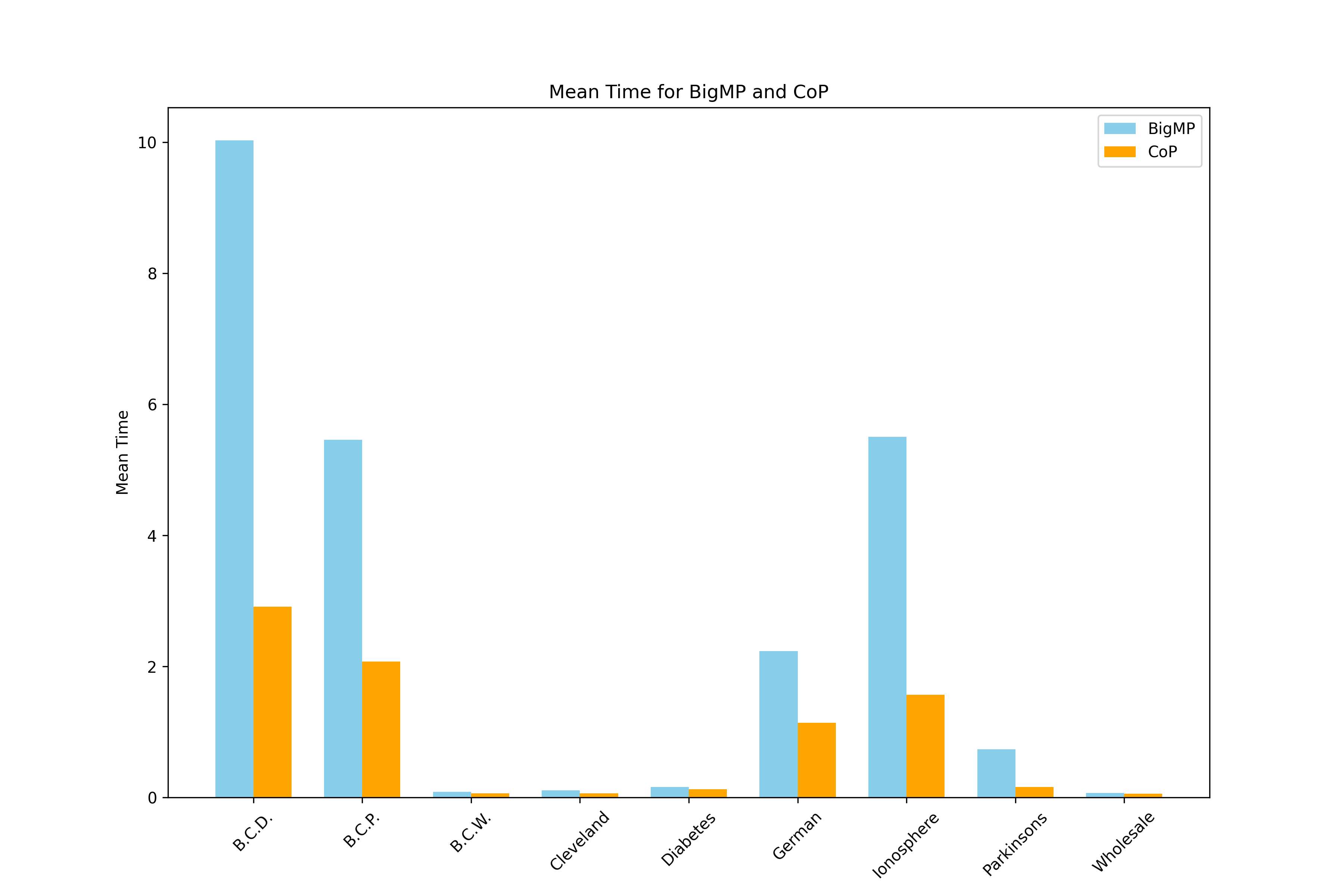}
        \caption{Comparison of solving BigMP and CoP with \texttt{Gurobi} for the \red{small} datasets: mean computational times (in seconds)}
    \label{fig:mean_time_small_data}
\end{figure}

\red{For the first set of results based on the small size datasets, we generated different instances ranging the integer $B$ from $1$ to $n$ and  $C$ values in the set $\{10^r : r \irg{-3}3\}$. 
We present the results  in Figure~\ref{fig:mean_time_small_data}, where we report for each dataset, the average computation times for solving the MIQP formulations, CoP and BigMP, on all the instances obtained varying $B,C$. Indeed, for this first set of problems, \texttt{Gurobi} was able to reach the global optimum in all the  cases, although the formulation CoP \eqref{prob:l2-FS-SVM2} was solved faster than BigMP~\eqref{prob:l2-FS-SVM1} by a factor of 1.4 on average across these datasets.}

\begin{table}[!ht]
\centering
\begin{tabular}{lccccc}

& & \multicolumn{2}{c}{BigMP} & \multicolumn{2}{c}{CoP} \\
    \hline
    \hline
Dataset & $B$ & ObjFun & MipGap  &  ObjFun & \ MipGap \\
\hline
\hline
\multirow{3}{*}{Arrhythmia} 
 & 10 & 2398.06 & \textbf{0.64} & \textbf{2151.56} & 0.82 \\
 & 20 & 2047.67 & \textbf{0.67} & \textbf{1712.64} & 0.81 \\
 & 30 & 1733.28 & \textbf{0.67} & \textbf{1441.63} & 0.78 \\
\hline
\multirow{3}{*}{Colorectal}
 & 10 & 2.99 & \textbf{0.64} & \textbf{2.12} & 0.98 \\
 & 20 & \textbf{0.64} & \textbf{0.13} & 0.68 & 0.93 \\
 & 30 & \textbf{0.41} & \textbf{0.04} & 0.43 & 0.89 \\
\hline
\multirow{3}{*}{DLBCL}
 & 10 & 0.63 & \textbf{1.00} & \textbf{0.56} & \textbf{1.00} \\
 & 20 & 0.18 & \textbf{0.99} & \textbf{0.15} & \textbf{0.99} \\
 & 30 & 0.16 & 0.99 & \textbf{0.09} & \textbf{0.98} \\
\hline
\multirow{3}{*}{Lymphoma} 
 & 10 & 18.07 & 1.00 & \textbf{1.16} & \textbf{0.99} \\
 & 20 & 0.71 & 0.99 & \textbf{0.29} & \textbf{0.98} \\
 & 30 & 0.34 & 0.98 & \textbf{0.16} & \textbf{0.96} \\
\hline
\multirow{3}{*}{{Madelon}} 
 & 10 & \textbf{16680.48} & \textbf{0.32} & 16800.07 & 0.33 \\
 & 20 & \textbf{16419.87} & \textbf{0.32} & 18558.98 & 0.40 \\
 & 30 & \textbf{16075.70} & \textbf{0.30} & 18425.44 & 0.39 \\
\hline
\multirow{3}{*}{{Mfeat}} 
 & 10 & 0.92 & \textbf{0.59} & \textbf{0.68} & 0.93 \\
 & 20 & \textbf{0.26} & \textbf{0.21} & \textbf{0.26} & 0.81 \\
 & 30 & \textbf{0.18} & \textbf{0.15} & \textbf{0.18} & 0.71 \\
\hline

\end{tabular}
\caption{\red{Comparison of solving BigMP and CoP with \texttt{Gurobi} for the large datasets: objective function (ObjFun)  values and optimal gap values (MipGap)}}
\label{tab: gurobi on big data} 
\end{table}


\red{For the results on datasets in Table~\ref{table:big_datasets}, we set the  hyperparameter $C=10$, and, similarly to~\cite{labbe2019mixed}, we ranged $B$ values in $\{10, 20, 30\}$.
In Table~\ref{tab: gurobi on big data} we show the computational performance of \texttt{Gurobi} and report the returned value of the objective function of problem \eqref{prob:l2-FS-SVM} (ObjFun) and the optimal gap value (MipGap) returned by Gurobi.}

 We can notice how the resulting optimization problems are much harder to solve in contrast to the problems related to the smaller datasets.  
 \red{Indeed for these results a} time limit of 3600 seconds was set, and the final MipGap values for the majority of these datasets are notably high. This is certainly due to the fact that the lower bound values computed in the Branch $\&$ Bound tree of \texttt{Gurobi} are not very effective for this type of formulations. This is particularly the case for model CoP whose final MipGap values tend to be worse than the ones of model BigMP. In contrast, \texttt{Gurobi} manages to find much better incumbent solutions when solving model CoP for all datasets, except for Colorectal and Madelon.

\subsection{Results on the heuristics}

In this section we analyse the performance of the heuristic algorithm we proposed and we compare it with~\texttt{Gurobi}. We have two versions of the heuristic algorithm depending on whether the Local Search or the Kernel Search upper bound strategies are used. We denote by \texttt{H-LS} the heuristic algorithm which implements the Local Search strategy, and by \texttt{H-KS} the heuristic implementing Kernel Search.

Table~\ref{table: h1 vs h2 large datasets} shows the optimization performance of \texttt{H-LS} and \texttt{H-KS}. Both heuristics \red{were} given a time limit of 600 seconds. In particular for \texttt{H-LS} we set the parameter $k=10$, which means that \red{ we solved the CoP model using  the $B+10$ most promising features based on the information of the optimal solution to the relaxation}. For \texttt{H-KS} we also set parameter $s=10$, which means that the overall set of features was divided in subsets of ten features each. A smaller time limit of 60 seconds was set for the resolution of each CoP subproblem solved in the Kernel Search strategy. 

As expected, \texttt{H-KS} manages to find better upper bounds compared to \texttt{H-LS}, at the cost of generally longer computation times. In some instances, such as Lymphoma, Madelon and Mfeat, both heuristics find the same solutions.  The time limit of 600 seconds was reached by \texttt{H-KS} for the instances related to datasets Madelon and Mfeat. For the rest of the datasets, the average computational times do not exceed 300 seconds for both heuristics. In general, compared to the objective function values of the solutions found by \texttt{Gurobi} in 1 hour~(Table \ref{tab: gurobi on big data}), the solutions found by the heuristics are much better.

\begin{table}[!ht]
\centering
\begin{tabular}{lcrrrr}
& & \multicolumn{2}{c}{\texttt{H-LS}} & \multicolumn{2}{c}{\texttt{H-KS}} \\
\hline
\hline
Dataset & $B$ & ObjFun & Time & ObjFun & Time \\
\hline
\hline
 \multirow{3}{*}{Arrhythmia} 
 & 10 & 2491.04 & \textbf{5.37} & \textbf{2142.67} & 41.24 \\
 & 20 & 1982.08 & \textbf{9.55} & \textbf{1652.83} & 168.83 \\
 & 30 & 1729.55 & \textbf{7.59} & \textbf{1380.19} & 311.62 \\
\hline
\multirow{3}{*}{Colorectal} 
 & 10 & 2.13 & \textbf{15.38} & \textbf{2.04} & 15.88 \\
 & 20 & 0.67 & \textbf{20.12} & \textbf{0.65} & 57.51 \\
 & 30 & 0.411 & \textbf{37.26} & \textbf{0.406} & 69.46 \\
\hline
\multirow{3}{*}{DLBCL} 
 & 10 & 0.25 & 51.01 & \textbf{0.23} & \textbf{26.56} \\
 & 20 & 0.104 & \textbf{46.97} & \textbf{0.101} & 67.74 \\
 & 30 & 0.065 & 134.95 & \textbf{0.063} & \textbf{117.29} \\
\hline
\multirow{3}{*}{Lymphoma}
 & 10 & \textbf{0.69} & \textbf{28.73} & \textbf{0.69} & 98.64 \\
 & 20 & 0.24 & \textbf{40.87} & \textbf{0.23} & 54.42 \\
 & 30 & 0.14 & \textbf{98.25} & \textbf{0.13} & 109.33 \\
 \hline
\multirow{3}{*}{Madelon} 
 & 10 & 16660.85 & \textbf{25.07} & \textbf{16615.43} & 607.32 \\
 & 20 & \textbf{16236.87} & \red{\textbf{48.57}} & \textbf{16236.87} & 612.62 \\
 & 30 & \textbf{15917.52} & \red{\textbf{176.92}} & 15917.99 & 605.93 \\
\hline
\multirow{3}{*}{Mfeat} 
 & 10 & 0.70 & \textbf{21.27} & \textbf{0.66}  & 341.87 \\
 & 20 & \textbf{0.25} & \textbf{76.68} & \textbf{0.25} & 621.16 \\
 & 30 & \textbf{0.17} & \textbf{276.90} & \textbf{0.17} & 606.96 \\
 \hline
\end{tabular}
\caption{Results on large datasets of the Heuristic algorithm implementing either the Local search (\texttt{H-LS}) or the Kernel Search (\texttt{H-KS}) as the upper bound strategies}
\label{table: h1 vs h2 large datasets}
\end{table}

Another set of results we carried out was the one of using \texttt{H-KS} as a warm start heuristic in order to evaluate if \texttt{Gurobi} was able to find better solution. In Table~\ref{table: heurs vs gurobi big data} we compare the optimal value of the solution of \texttt{Gurobi} after 1 hour for model CoP, the value of the \texttt{H-KS} heuristic, and the value of the solution found by \texttt{Gurobi} using as a warm start solution the one of  \texttt{H-KS} (denoted by  CoP$^*$). 

Regarding Table~\ref{table: heurs vs gurobi big data}, we can see that \texttt{Gurobi} was unable to find a better solution than the heuristic after 1 hour for all datasets, apart from Arrhythmia, and this is probably due to the fact that the heuristic solutions are already optimal or nearly-optimal. In this table, column $d_f$ indicate the percentage of improvement of the objective function found by \texttt{Gurobi} using \texttt{H-KS} starting solution, against the objective function of the solution found by \texttt{Gurobi} alone.

\begin{table}[!ht]
\centering
\begin{tabular}{lcrrrc}
& & \multicolumn{3}{c}{ObjFun} & \\
\hline
\hline
Dataset & $B$ & CoP & \texttt{H-KS} & CoP$^*$ & $d_f$ \\
\hline
\hline
 \multirow{3}{*}{Arrhythmia} 
 & 10 & 2151.56 & 2142.67 & \textbf{2098.55} & 2\% \\
 & 20 & 1712.64 & \textbf{1652.83} & \textbf{1652.83} & 3\% \\
 & 30 & 1441.63 & 1380.19 & \textbf{1366.94} & 5\% \\
\hline
\multirow{3}{*}{Colorectal}
 & 10 & 2.12 & \textbf{2.04} & \textbf{2.04} & 4\% \\
 & 20 & 0.68 & \textbf{0.65} & \textbf{0.65} & 5\% \\
 & 30 & 0.43 & \textbf{0.41} & \textbf{0.41} & 6\% \\
\hline
\multirow{3}{*}{DLBCL} 
 & 10 & 0.56 & \textbf{0.23} & \textbf{0.23} & 60\% \\
 & 20 & 0.15 & \textbf{0.10} & \textbf{0.10} & 35\% \\
 & 30 & 0.09 & \textbf{0.06} & \textbf{0.06} & 29\% \\
\hline
\multirow{3}{*}{Lymphoma} 
 & 10 & 1.16 & \textbf{0.69} & \textbf{0.69} & 41\% \\
 & 20 & 0.29 & \textbf{0.22} & \textbf{0.22} & 21\% \\
 & 30 & 0.16 & \textbf{0.13} & \textbf{0.13} & 16\% \\
 \hline
\multirow{3}{*}{Madelon}
 & 10 & 16800.07 & \textbf{16615.43} & \textbf{16615.43} & 1\% \\
 & 20 & 18558.98 & \textbf{16236.86} & \textbf{16236.86} & 13\% \\
 & 30 & 18425.44 & \textbf{15917.51} & \textbf{15917.51} & 14\% \\
\hline
\multirow{3}{*}{Mfeat} 
 & 10 & 0.68 & \textbf{0.66} & \textbf{0.66} & 3\% \\
 & 20 & 0.26 & \textbf{0.25} & \textbf{0.25} & 4\% \\
 & 30 & \red{\textbf{0.17}} & \textbf{0.17} & \textbf{0.17} & 3\% \\
 \hline
\end{tabular}
\caption{Comparison of the objective functions of CoP returned by \texttt{Gurobi} with and without using the Heuristic Kernel Search algorithm as a warm start for large datasets}
\label{table: heurs vs gurobi big data}
\end{table}

\begin{figure}[h!]
        \centering
        \includegraphics[width=16cm]{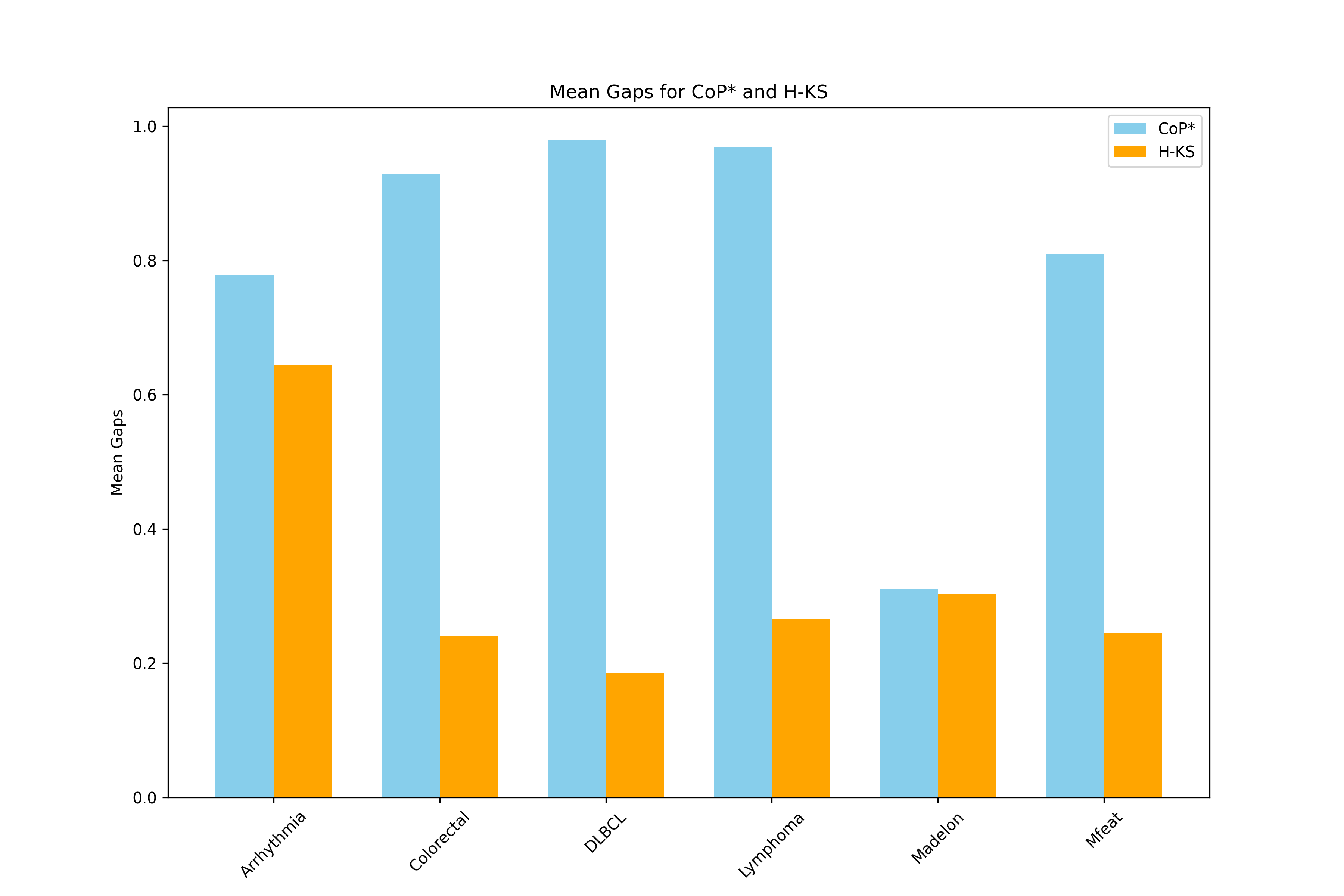}
        \caption{Comparison of the mean MipGap values of the solution found by  \texttt{Gurobi} solving CoP after 1 hour, and the Relaxed Gap of the Heuristic Kernel Search solution using the solution of DSCoP as a lower bound}
    \label{fig:gaps comparison}
\end{figure}

\begin{table}[!ht]
\centering
\begin{tabular}{lcrrr}
& & \multicolumn{3}{c}{Initial Bounds}  \\
\hline
\hline
Dataset & $B$ & \text{DSCoP} &BigMP & CoP \\
\hline
\hline
\multirow{3}{*}{Arrhythmia} 
 & 10 & \textbf{754.46} & 451.10 & 274.79 \\
 & 20 & \textbf{622.44} & 395.92 & 274.79 \\
 & 30 & \textbf{538.86} & 303.66 & 274.79 \\
\hline
\multirow{3}{*}{Colorectal} 
 & 10 & \textbf{1.07} & 0.05 & 0.05 \\
 & 20 & \textbf{0.54} & 0.05 & 0.05 \\
 & 30 & \textbf{0.38} & 0.05 & 0.05 \\
\hline
\multirow{3}{*}{DLBCL} 
 & 10 & \textbf{0.16} & 0.00 & 0.00 \\
 & 20 & \textbf{0.08} & 0.00 & 0.00 \\
 & 30 & \textbf{0.06} & 0.00 & 0.00 \\
\hline
\multirow{3}{*}{Lymphoma} 
 & 10 & \textbf{0.35} & 0.00 & 0.01 \\
 & 20 & \textbf{0.18} & 0.00 & 0.01 \\
 & 30 & \textbf{0.12} & 0.00 & 0.01 \\
\hline
\multirow{3}{*}{Madelon} 
 & 10 & \textbf{11243.14} & 0.00 & 0.00 \\
 & 20 & \textbf{11209.44} & 0.00 & 0.00 \\
 & 30 & \textbf{11197.73} & 0.00 & 0.00 \\
\hline
\multirow{3}{*}{Mfeat} 
 & 10 & \textbf{0.36} & 0.05 & 0.05 \\
 & 20 & \textbf{0.20} & 0.05 & 0.05 \\
 & 30 & \textbf{0.15} & 0.05 & 0.05 \\
\hline
\end{tabular}
\caption{Lower bounds comparison between the SDP one and the starting bound found by \texttt{Gurobi} for the BigMP and CoP models for different datasets}
\label{table: lower_bounds_comparison}
\end{table}

\subsection{Results on the exact approach}

We now present the results on the Exact procedure. As shown in Table \ref{tab: gurobi vs sfs-e}, the  exact algorithm performs much better than $\texttt{Gurobi}$ with a warm start heuristic,  both on optimal value and computation time. In particular, the exact algorithm is able to globally solve four problems under a time limit of 3600 seconds. \red{It is worth mentioning that the initial lower bounds obtained through the relaxation \eqref{prob:dsdp-l2-FS-SVM2} for these problems are already quite strong. In the process of the exact procedure, after several updates of the set $K$, the semi-relaxation SR-DLMP($K$) provides an almost tight bound, allowing the algorithm to terminate within the time limit.} For the unsolved problems, the exact algorithm provides a better objective value and a smaller MipGap.  
\begin{table}[!ht]
\centering
\begin{tabular}{lcrccrcc}
& & \multicolumn{3}{c}{CoP$^*$} & \multicolumn{3}{c}{\texttt{Exact KS}} \\
\hline
\hline
Dataset & $B$ & ObjFun & Time & MipGap &  ObjFun & Time &  MipGap  \\
\hline
\hline

\multirow{3}{*}{Arrhythmia}
& 10 & 2098.55 & 3600 & 0.82 & \textbf{2003.98} & 3600 & \textbf{0.22}  \\
& 20 & 1652.83  & 3600 & 0.80 & \textbf{1632.32} & 3600 & \textbf{0.27} \\
& 30 & 1366.94 & 3600 & 0.77 & \textbf{1362.75} & 3600 & \textbf{0.32}  \\
\hline
\multirow{3}{*}{Colorectal} 
& 10 & \red{\textbf{2.04}} & 3600 & 0.98 & \textbf{2.04}  & \textbf{407.2} & \textbf{0}   \\
& 20 & \red{\textbf{0.65}} & 3600 & 0.93 & \textbf{0.65} & \textbf{431.5} & \textbf{0}   \\
& 30 & \red{\textbf{0.41}} & 3600 & 0.88 & \textbf{0.41}  & \textbf{407.6} & \textbf{0}   \\
 \hline
\multirow{3}{*}{DLBCL}  
& 10 &  \textbf{0.23} & 3600 & 0.99 & \textbf{0.23} & \textbf{451.3} &  \textbf{0} \\
& 20 &  \textbf{0.10} & 3600 & 0.98 & \textbf{0.10} & \textbf{542.7} &  \textbf{0} \\
& 30 &  \textbf{0.06} & 3600 & 0.97 & \textbf{0.06} & \textbf{746.2} &  \textbf{0} \\
\hline
\multirow{3}{*}{Lymphoma} 
& 10 & \textbf{0.69} & 3600  & 0.99 & \textbf{0.69}  & \textbf{314.6} & \textbf{0} \\
& 20 & \textbf{0.22} & 3600  & 0.97 & \textbf{0.22}  & \textbf{532.3} & \textbf{0} \\
& 30 & \textbf{0.13} & 3600  & 0.95 & \textbf{0.13}  & \textbf{823.4} & \textbf{0} \\
\hline
\multirow{3}{*}{{Madelon}} 
& 10 & 16615.43 & 3600 & 0.33  & \textbf{16499.38} & 3600 & \textbf{0.17}\\
& 20 & 16236.86 & 3600 & 0.31  & \textbf{16142.54} & 3600 & \textbf{0.18}\\
& 30 & 15917.51 & 3600 & 0.30 & \textbf{15765.83} & 3600 & \textbf{0.21} \\
\hline
\multirow{3}{*}{{Mfeat}}
& 10 & 0.66 & 3600 & 0.93 &  \textbf{0.64}   &  \textbf{397.5} & \textbf{0} \\
        & 20 & 0.25 & 3600 & 0.80 & \textbf{0.24}   & \textbf{573.5} & \textbf{0}  \\
& 30 & \textbf{0.17} & 3600 & 0.71 & \textbf{0.17}   & \textbf{1045.5} & \textbf{0}  \\
\hline
\end{tabular}
\caption{Comparison of the \texttt{Gurobi}  solution to the CoP model using the Heuristic Kernel Search solution as a warm start heuristic, and the Exact algorithm using the Kernel Search procedure as the upper bound strategy.}
\label{tab: gurobi vs sfs-e}
\end{table}

\subsection{Validation}

Finally, we conducted a series of experiments to compare our feature selection model \eqref{prob:l2-FS-SVM} which is based on the $\ell_2$-norm, \red{against other machine learning models in the literature. Firstly, we compare it against}  the feature selection model studied in \cite{labbe2019mixed} which employs the $\ell_1$-norm instead; \red{secondly we compare it with a model which adds an $\ell_1$ penalty term in the objective function (similarly to the elastic-net models proposed in \cite{wang2006doubly,zou2005regularization}). This approach, aimed at achieving solution sparsity, requires tuning of the penalty parameter. The aim of the first set of experiments is to validate, from a machine learning perspective, our model compared to the one by \cite{labbe2019mixed} which was already proven to be more effective than other SVM-based feature selection models. The second set of experiments is instead aimed at emphasizing the main difference between our model, which is based on a strict cardinality constraint, and an elastic-net-type model which has no direct control over the number of features that will be selected}

\red{As for the comparison with the $\ell_1$ model,} we varied the parameter $C$ among values in $\{10^r : r \irg{-3}3\}$. Regarding the parameter $B$, similarly to \cite{labbe2019mixed}, we ranged it among all possible values in-between  $1$ and $n$ for the small datasets, while, for the larger ones, $B$ varied in the set $\{10,20,30,40,50\}$. 

\begin{sidewaysfigure}[ht]
  \centering
  \begin{subfigure}[b]{0.3\textwidth}
    \includegraphics[width=\textwidth]{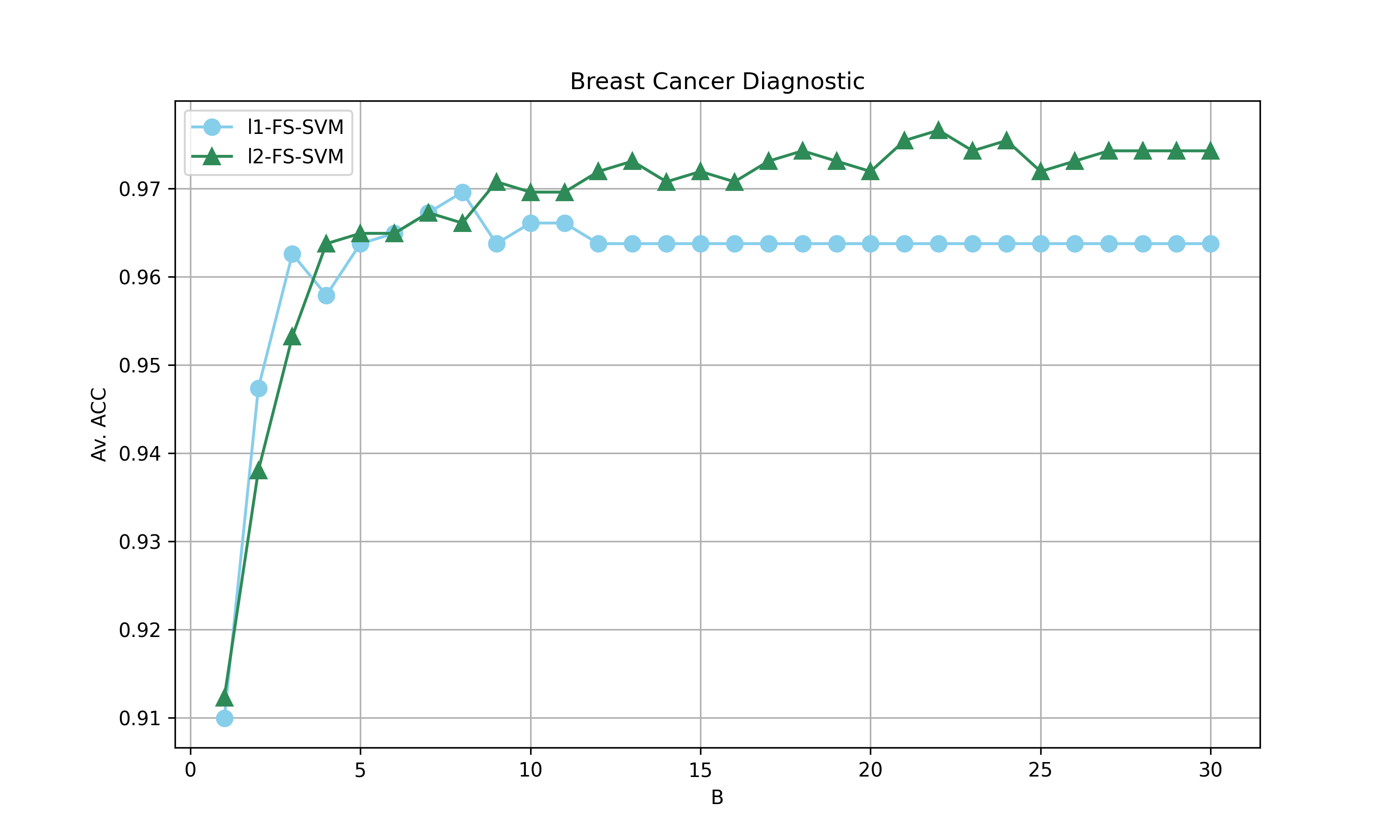}
  \end{subfigure}
  \hfill
  \begin{subfigure}[b]{0.3\textwidth}
    \includegraphics[width=\textwidth]{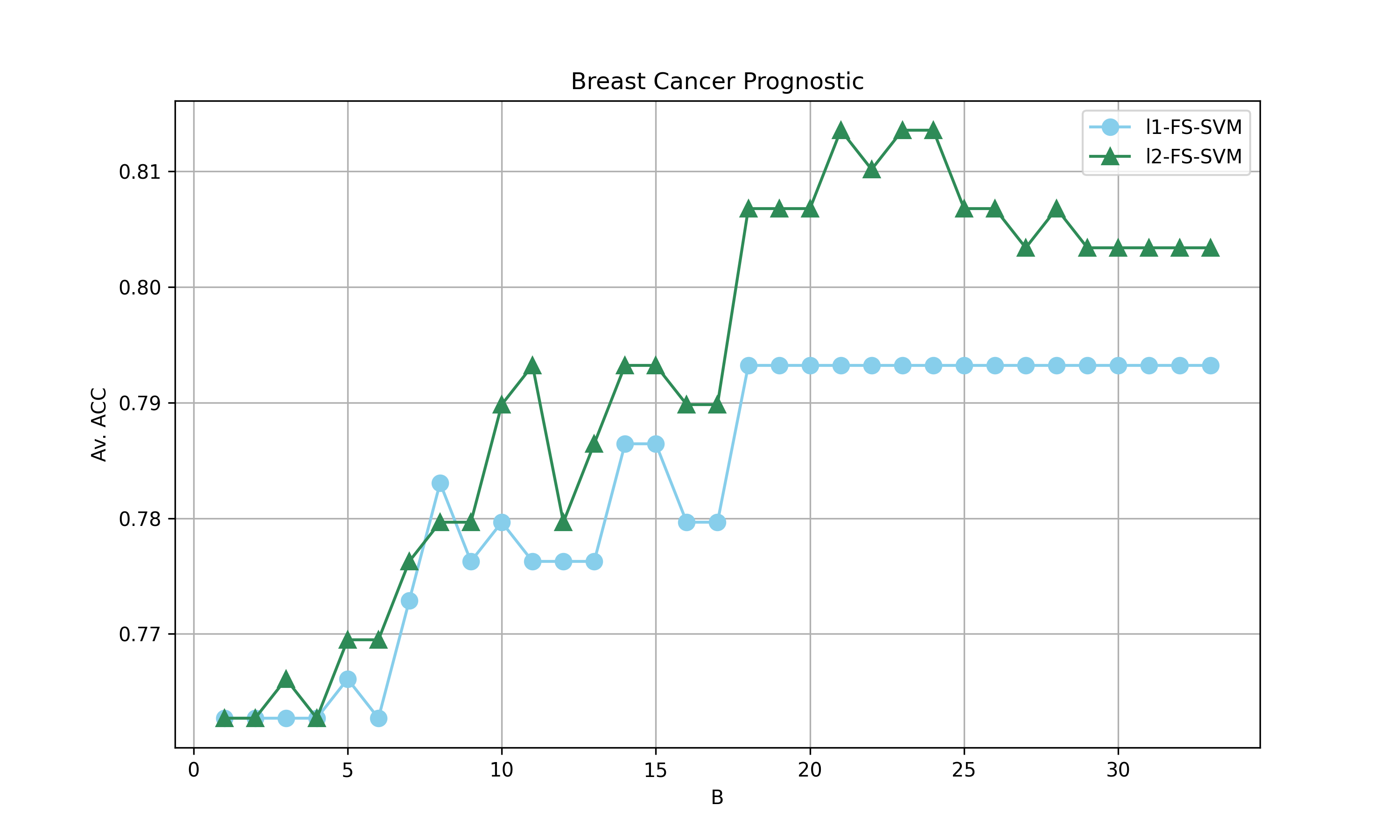}
  \end{subfigure}
  \hfill
  \begin{subfigure}[b]{0.3\textwidth}
    \includegraphics[width=\textwidth]{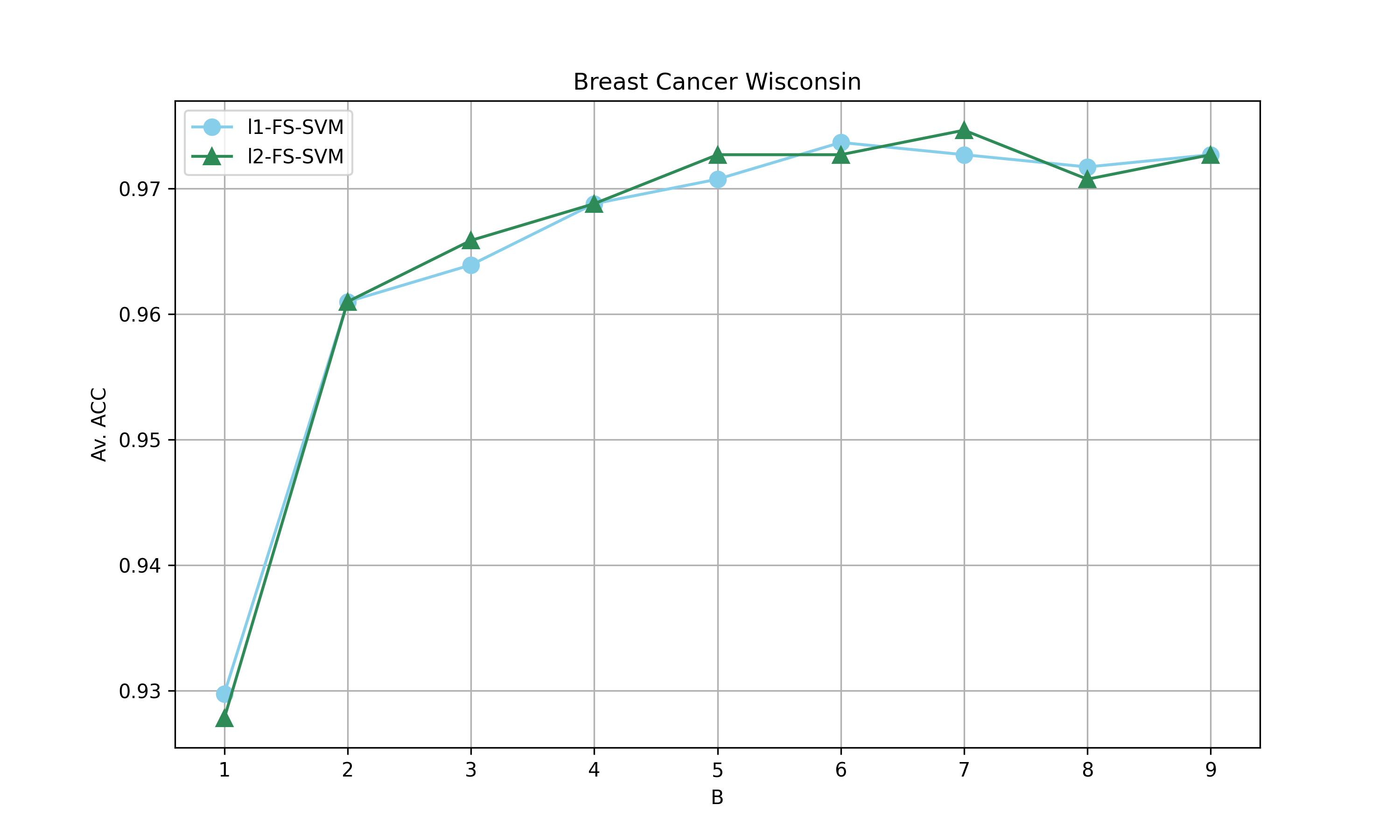}
  \end{subfigure}
  
  \begin{subfigure}[b]{0.3\textwidth}
    \includegraphics[width=\textwidth]{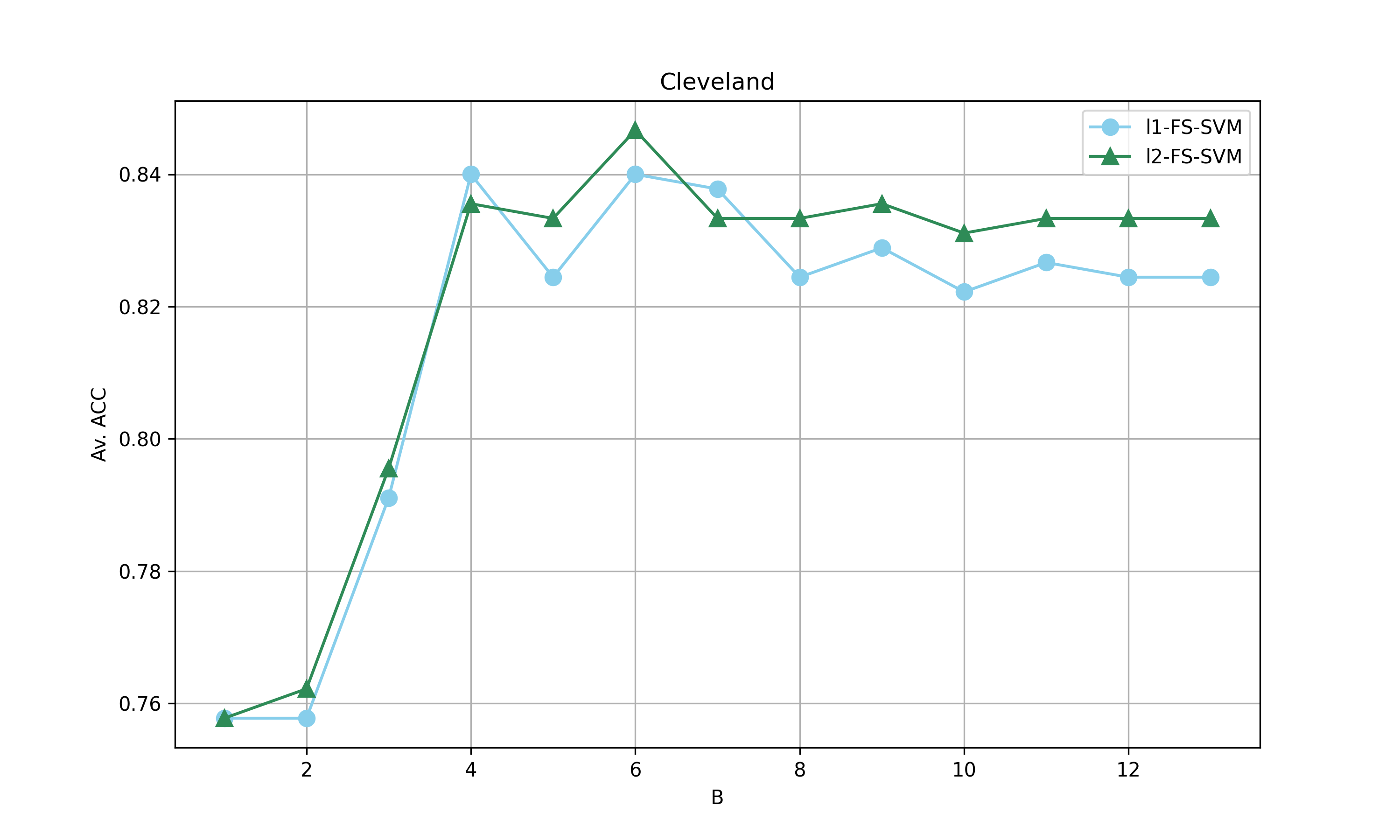}
  \end{subfigure}
  \hfill
  \begin{subfigure}[b]{0.3\textwidth}
    \includegraphics[width=\textwidth]{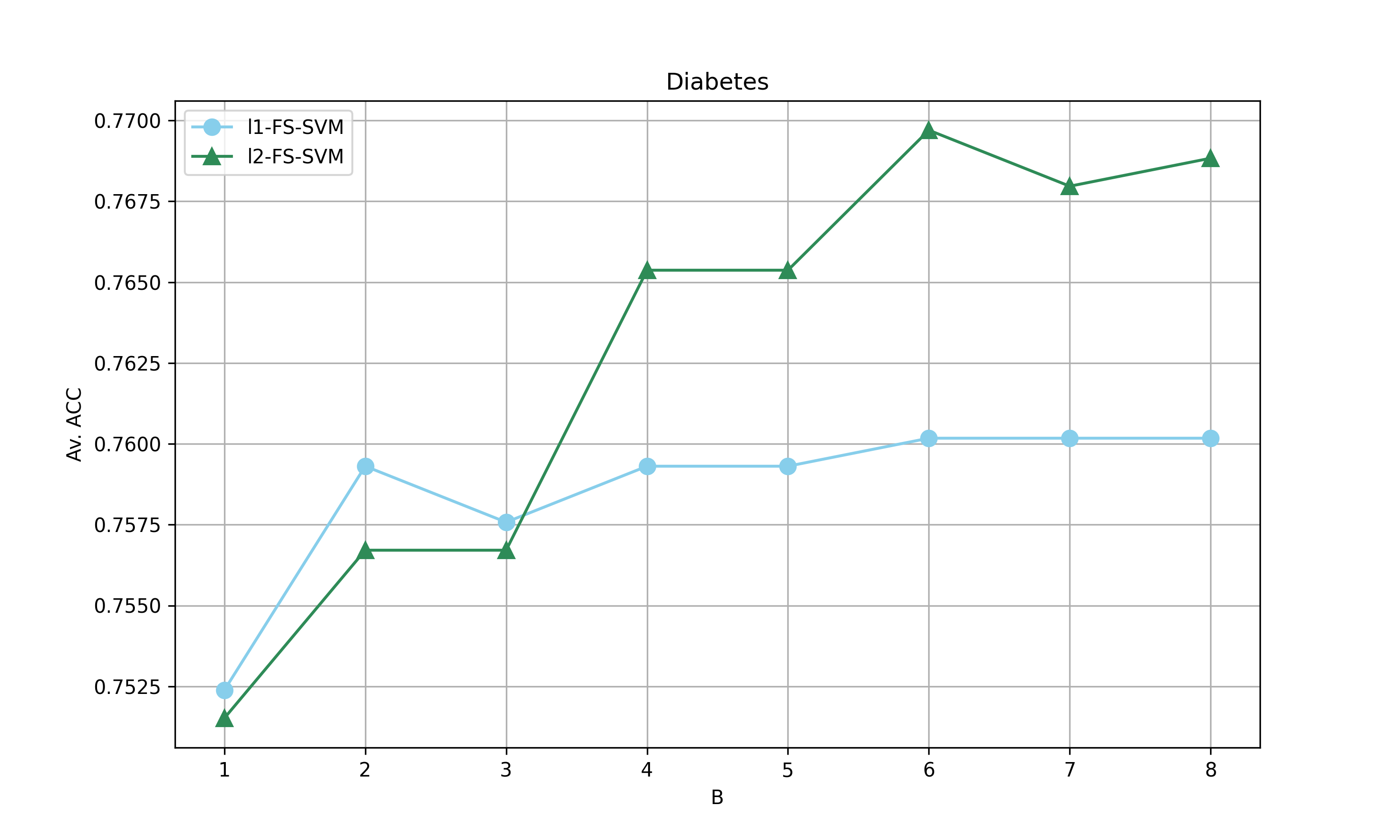}
  \end{subfigure}
  \hfill
  \begin{subfigure}[b]{0.3\textwidth}
    \includegraphics[width=\textwidth]{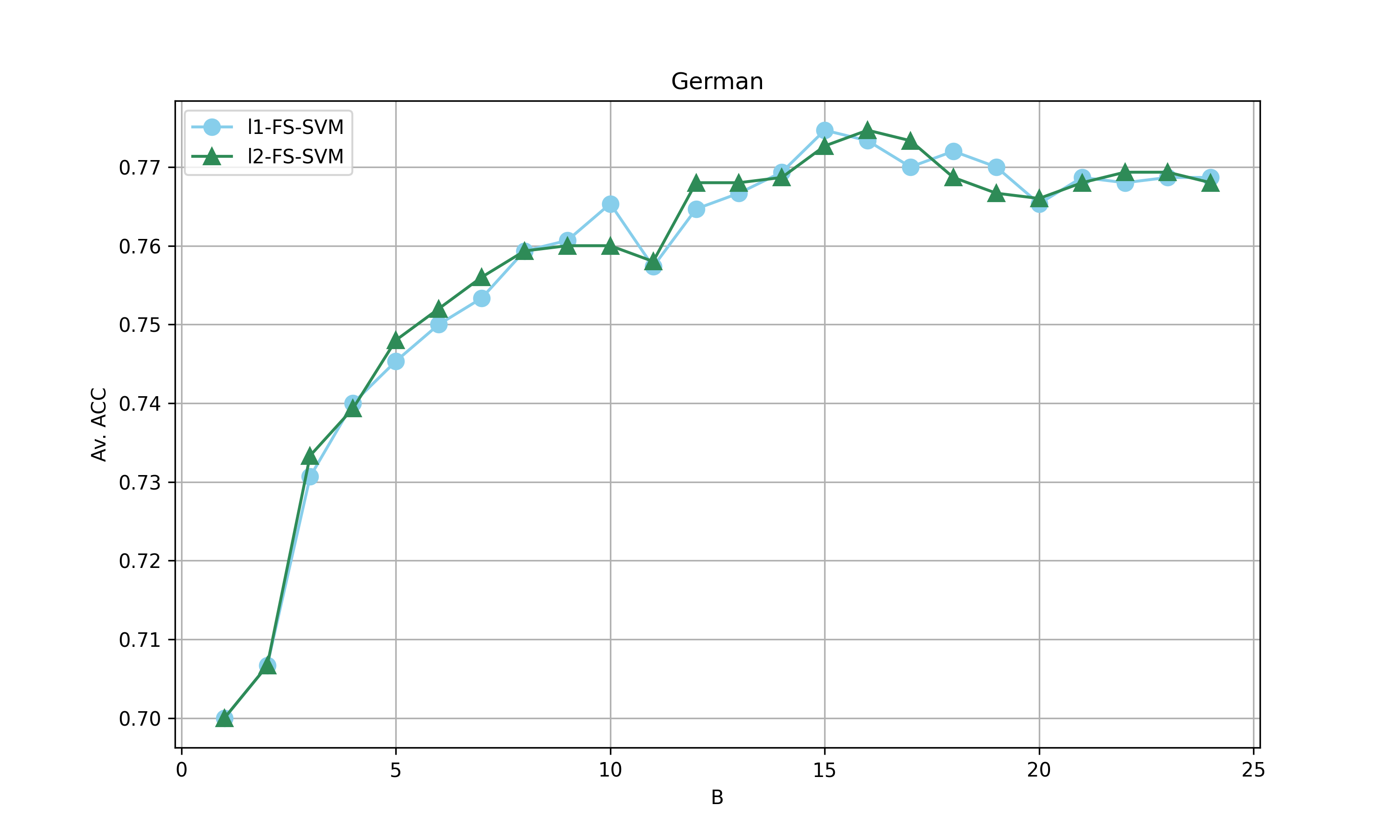}
  \end{subfigure}
  
  \begin{subfigure}[b]{0.3\textwidth}
    \includegraphics[width=\textwidth]{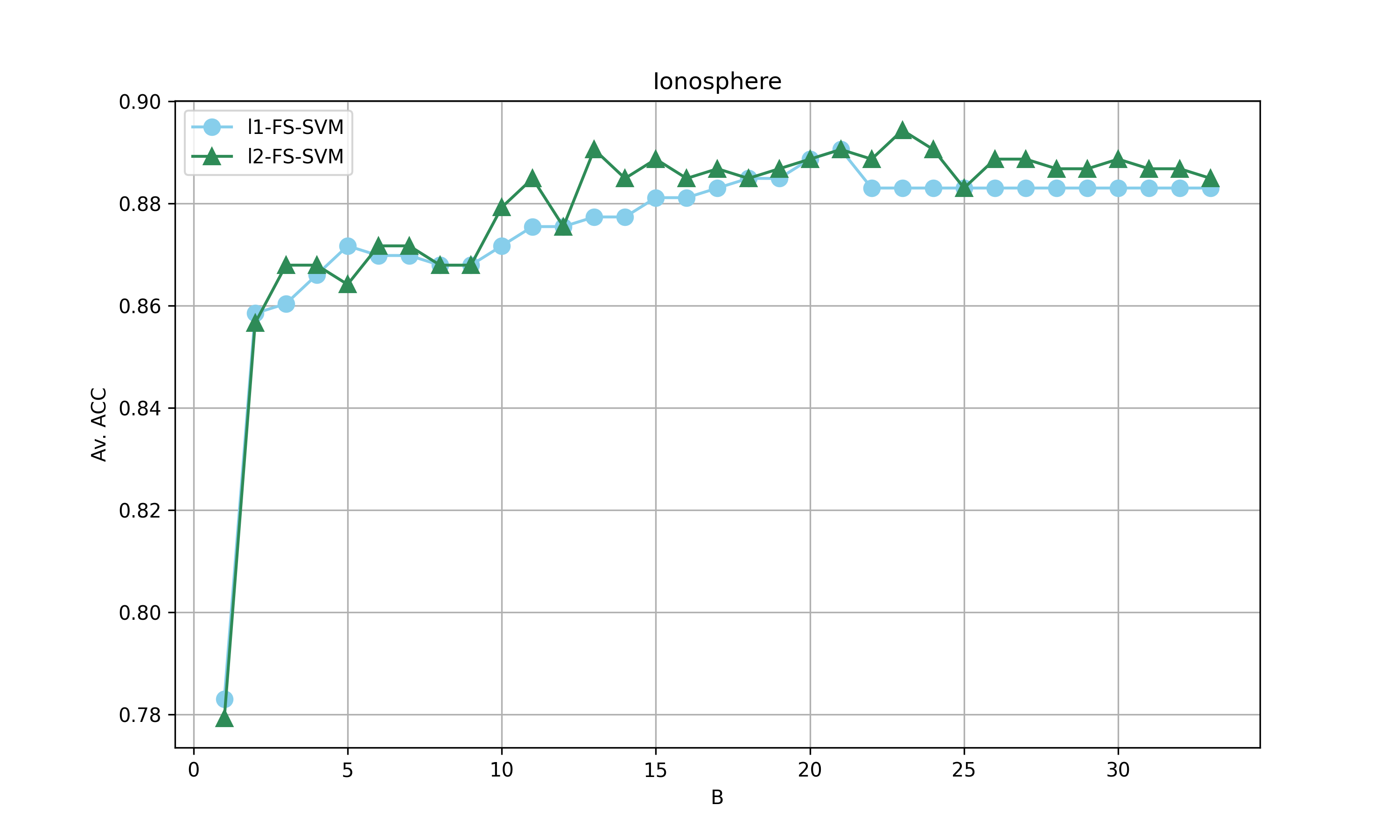}
  \end{subfigure}
  \hfill
  \begin{subfigure}[b]{0.3\textwidth}
    \includegraphics[width=\textwidth]{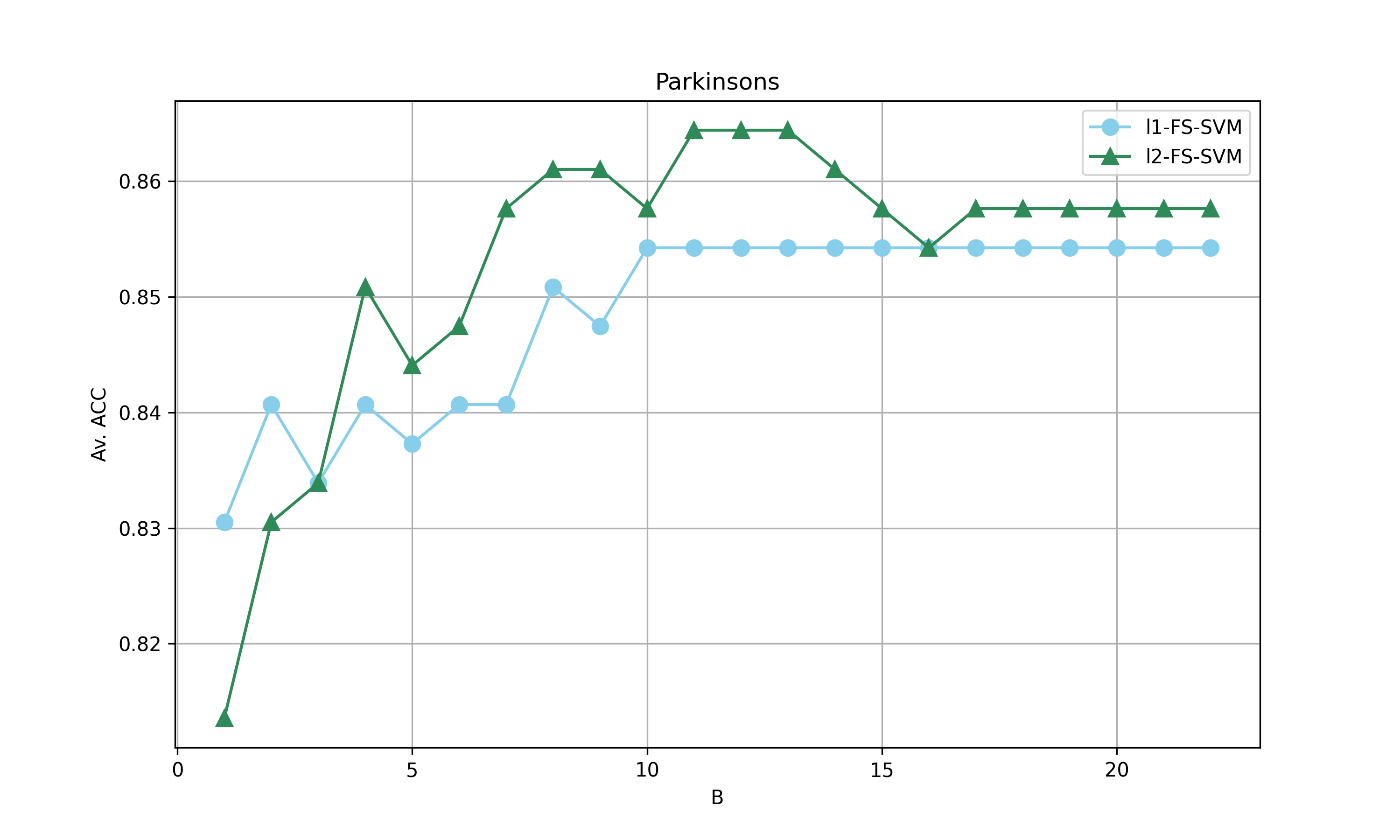}
  \end{subfigure}
  \hfill
  \begin{subfigure}[b]{0.3\textwidth}
    \includegraphics[width=\textwidth]{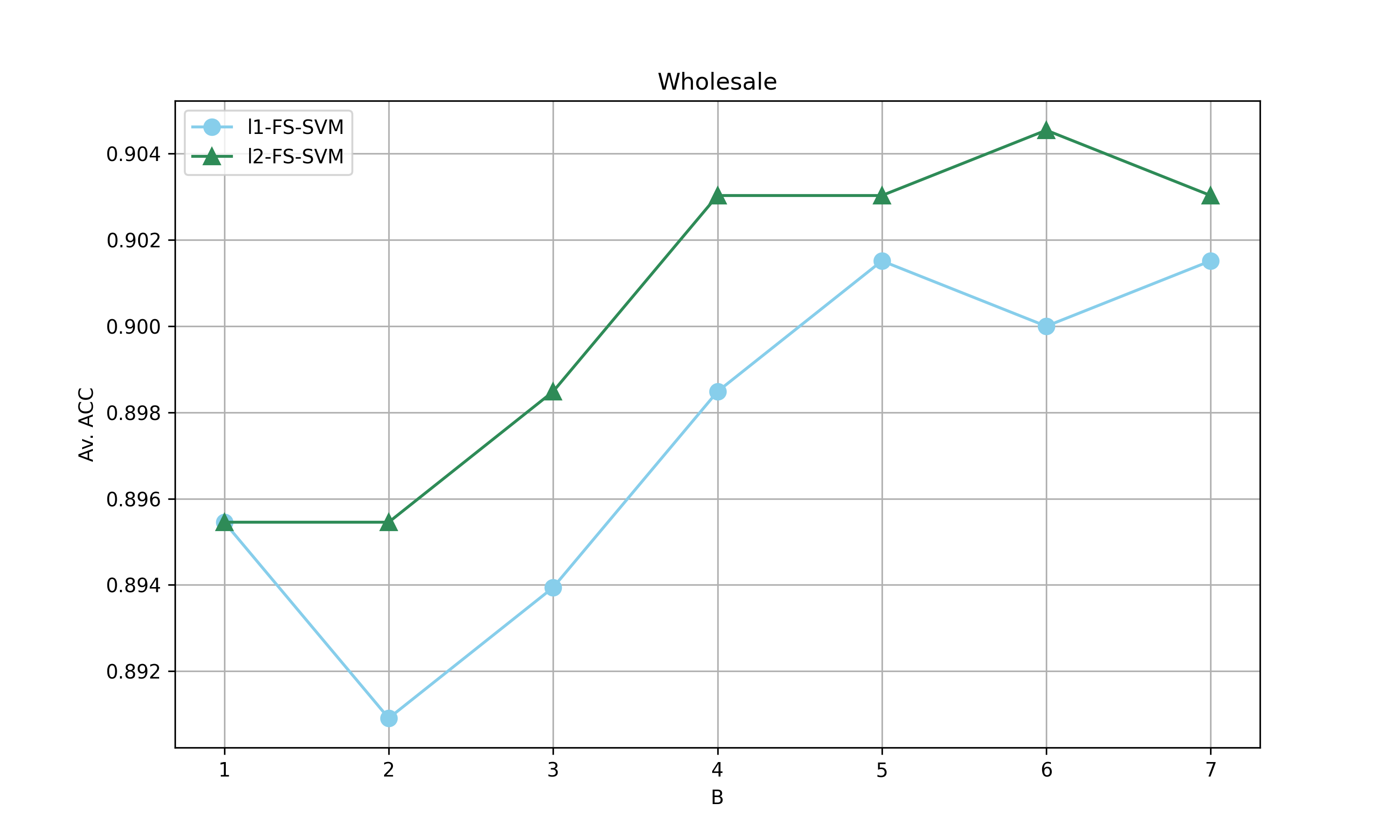}
  \end{subfigure}
  
  \caption{The $\ell_2$-FS-SVM (dark green) versus $\ell_1$-FS-SVM (light blue): best average validation ACC values over a 10-fold cross validation as $B$ changes - small datasets (problems solved to optimality).}
  \label{fig:AvACC small data}
\end{sidewaysfigure}

Under a ten-fold cross-validation framework, for each of the ten folds and for each $C$ and $B$ in the grid, we computed the accuracy (ACC) values on the portion of the dataset left for validation. We then computed the average ACC values among the \red{ten} folds. In Figure~\ref{fig:AvACC small data}, we plot the best mean ACC values for each $B$ in the grid for the small datasets. We recall that in this case, all the optimization problems could be solved to optimality by \texttt{Gurobi}. As we can see, in most of the cases, employing the $\ell_2$-norm leads to better generalization performance, resulting in higher average accuracy values.

\begin{sidewaysfigure}[ht]
  \centering
  \begin{subfigure}[b]{0.32\textwidth}
    \includegraphics[width=\textwidth]{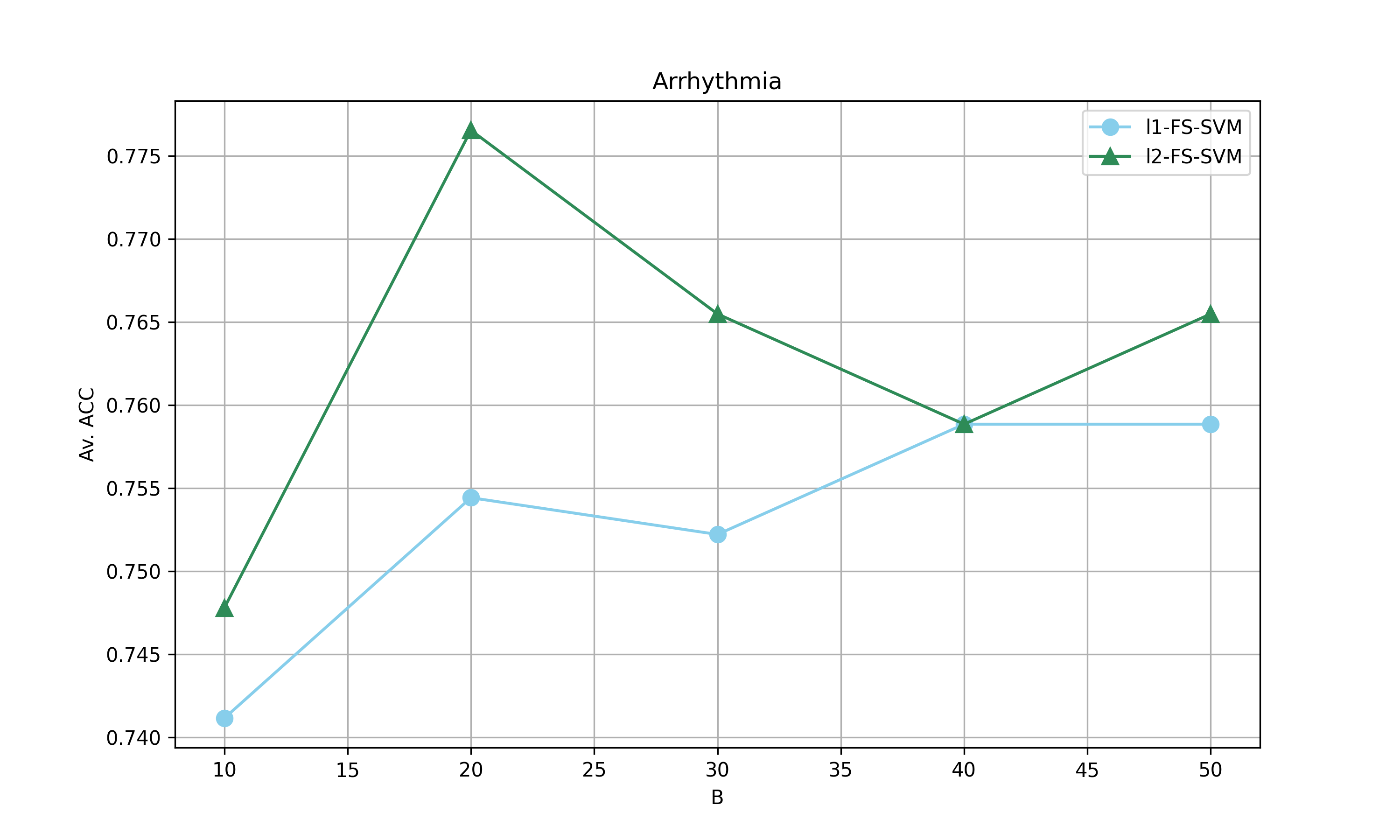}
  \end{subfigure}
  \hfill
  \begin{subfigure}[b]{0.32\textwidth}
    \includegraphics[width=\textwidth]{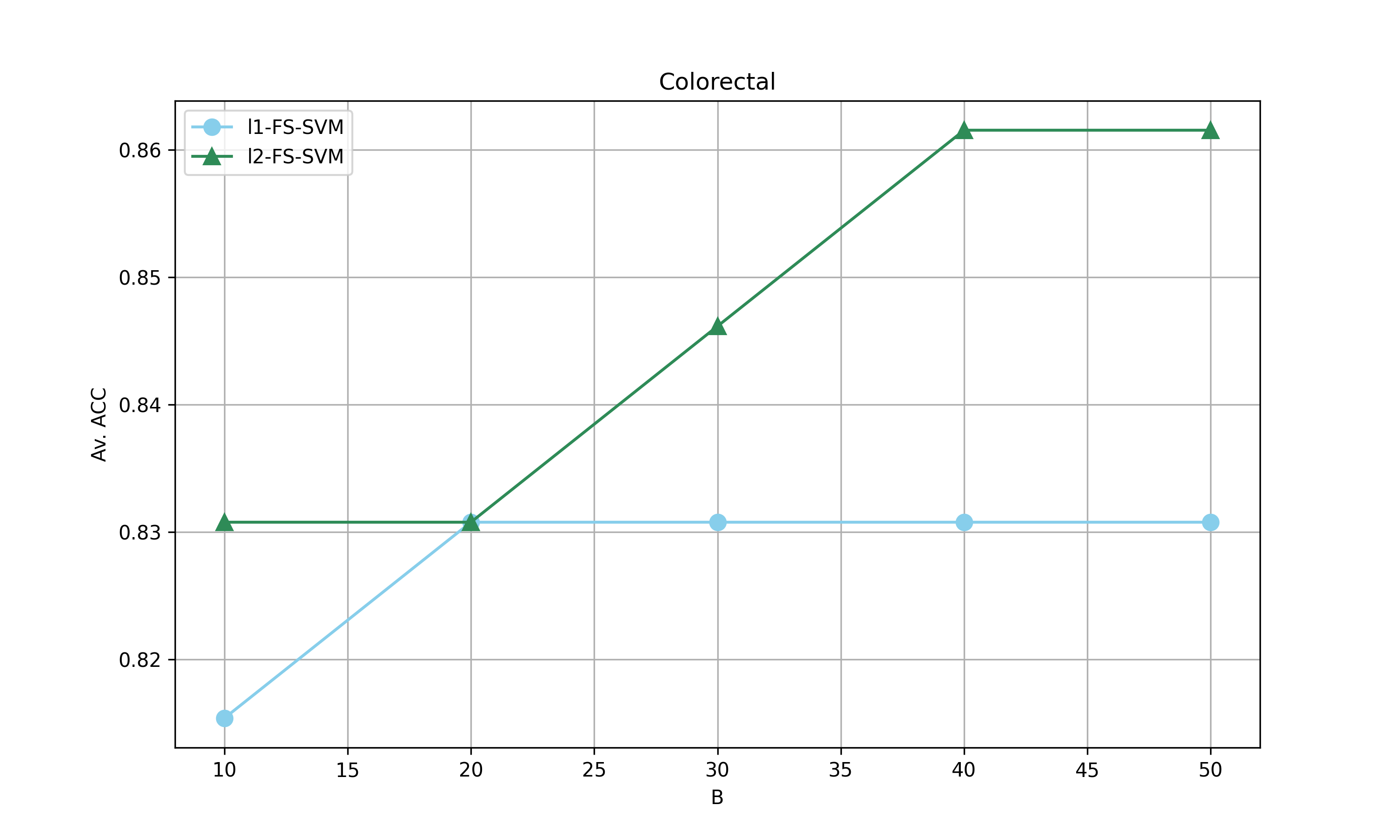}
  \end{subfigure}
  \hfill
  \begin{subfigure}[b]{0.32\textwidth}
    \includegraphics[width=\textwidth]{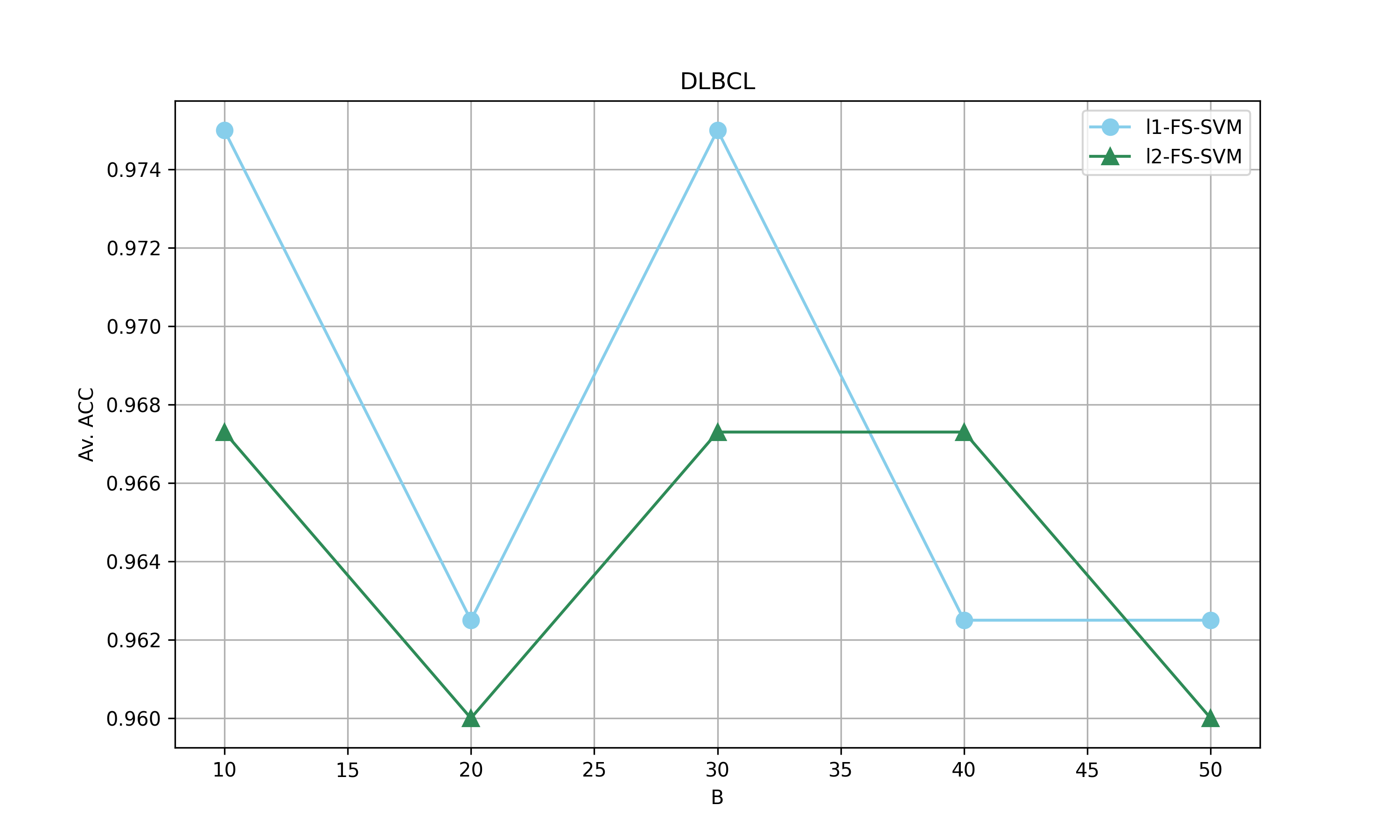}
  \end{subfigure}

  \begin{subfigure}[b]{0.32\textwidth}
    \includegraphics[width=\textwidth]{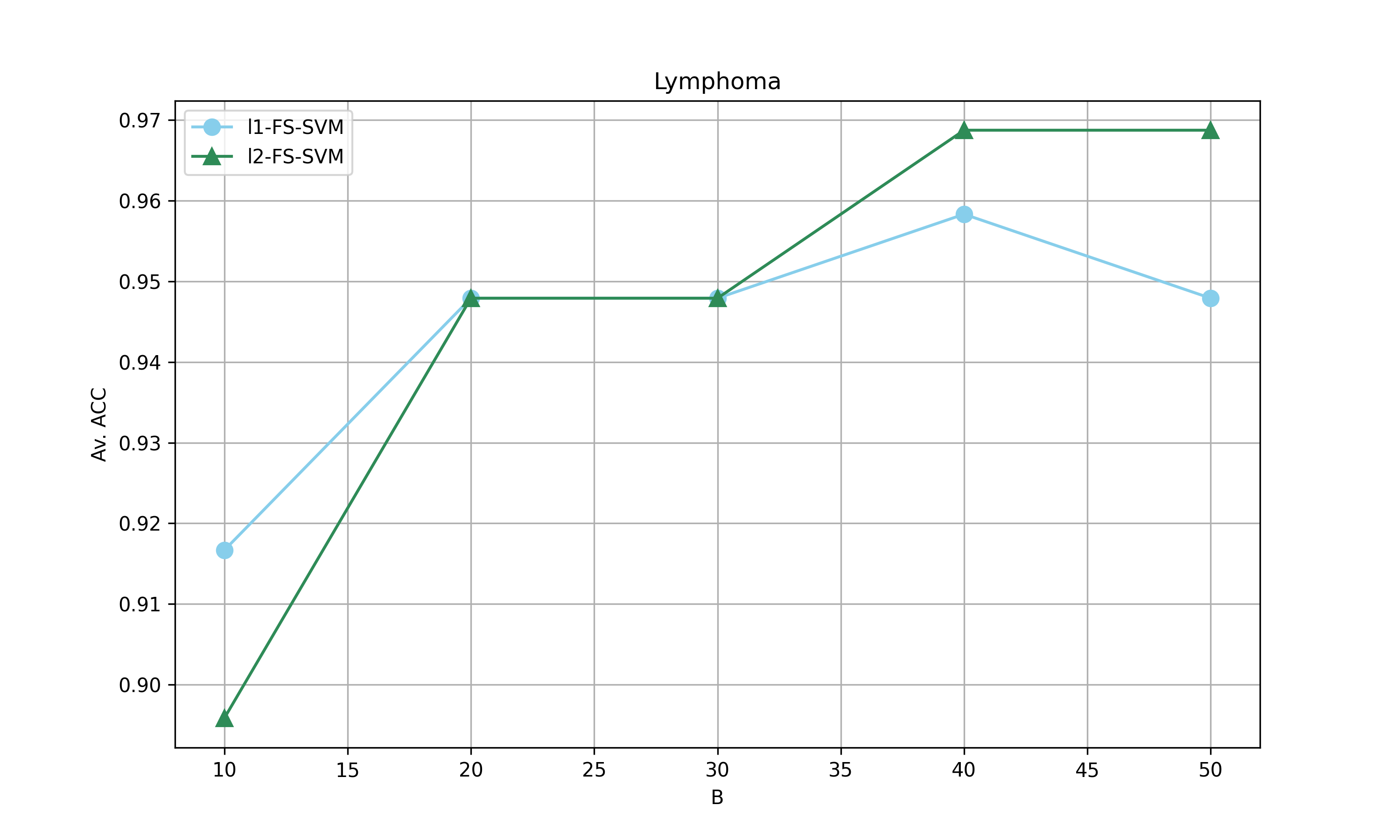}
  \end{subfigure}
  \hfill
  \begin{subfigure}[b]{0.32\textwidth}
    \includegraphics[width=\textwidth]{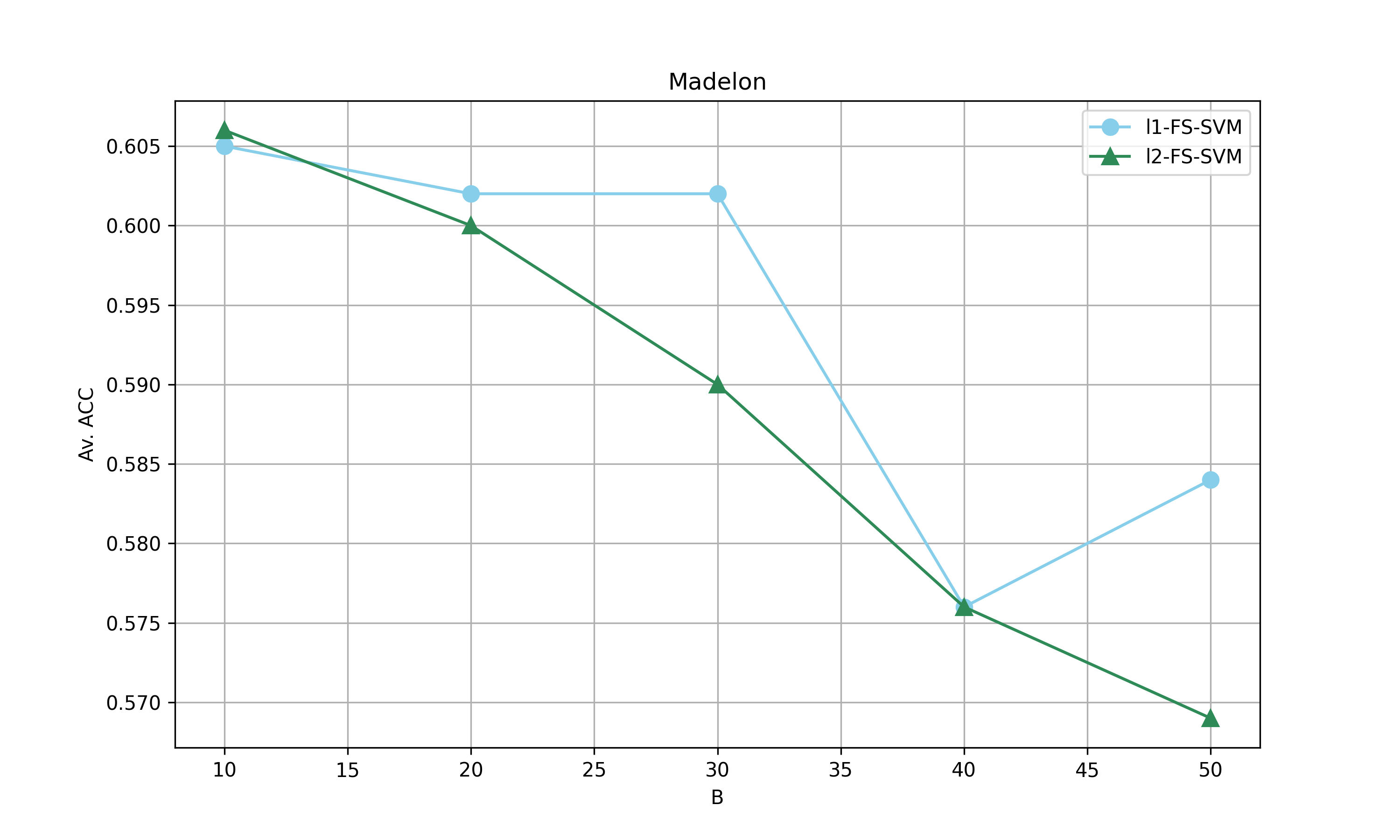}
  \end{subfigure}
  \hfill
  \begin{subfigure}[b]{0.32\textwidth}
    \includegraphics[width=\textwidth]{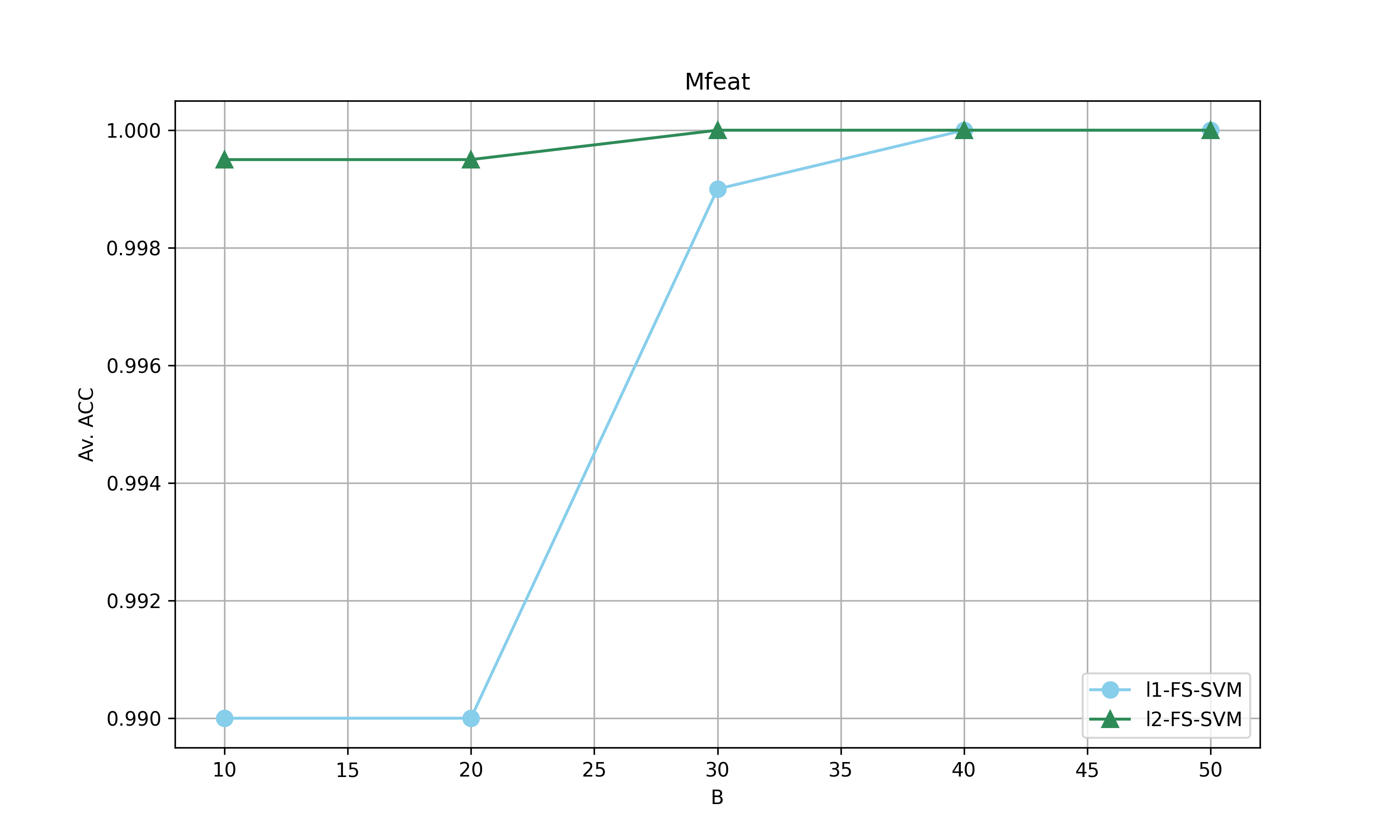}
  \end{subfigure}

  \caption{The $\ell_2$-FS-SVM (dark green) versus $\ell_1$-FS-SVM (light blue): best average validation ACC values over a 10-fold cross validation as $B$ changes  - large datasets (solutions returned by the Heuristic Kernel Search  with a time limit of 600 sec)}
  \label{fig:AvACC large data}
\end{sidewaysfigure}

Similar results for the large datasets are reported in Figure~\ref{fig:AvACC large data}. In this case, each optimization problem was solved with the Kernel Search strategy alone. This is because, as shown in the computational experience, the strategy was able to find in most cases very good, if not optimal, solutions. The same was done in~\cite{labbe2019mixed}. To differentiate the two models, the FS-SVM model we studied in this paper is here named as the $\ell_2$-FS-SVM model, while the one presented in \cite{labbe2019mixed} is referred to as the $\ell_1$-FS-SVM. The Kernel Search heuristic we implemented for the $\ell_1$-FS-SVM is the same as described in the relative paper. In order to conduct a fair comparison, we set the same time limit to 600 seconds for each procedure.  

\red{
As regards the comparison with the elastic-net-type SVM model \cite{wang2006doubly,zou2005regularization}, we perform experiments to check whether the approaches can obtain the same performance. For this set of results we did not carry out a $k$-fold cross validation but we partitioned the datasets in a training set and a test set ($80\%$ and $20\%$ of the datapoints respectively). We fixed the parameter $C$ to 10, and the budget value was set to $B=5$ for the set of the small datasets (Table \ref{table:small_datasets}), and to $B=10$ for the larger datasets (Table \ref{table:big_datasets}).
The results are reported in Tables~\ref{tab:E-NET1} and \ref{tab:E-NET2}.} 

\red{Our model, denoted as FS-SVM, always selects the desired number of features, due to the strict cardinality constraint. The elastic-net-type SVM model, denoted as EN-SVM, does not use the cardinality constraint, instead an additional penalization term $\alpha \|\w\|_1$ is added to the objective function. We ranged values of $\alpha$ in $\{ 2^i : i\irg 0 {10}\}$. For each dataset we reported the training and test accuracy values, as well as the number of features selected  (underlined values indicate that the number of selected features exceeds the limit $B$). If there are $\alpha$ values for the EN-SVM model which yielded a number of features not exceeding $B$, we chose the model with the best training accuracy, otherwise we chose the model with the lowest number of features selected. In the latter case, increasing $\alpha$ values for the selected $\alpha$, yielded models in which no features were selected (i.e. for those $\alpha$ values, the optimal $\w$ had all zero components).}

\red{The results show that the FS-SVM models consistently achieve better ACC values, not just on the training set, but also on the test set. In most of the experiments, the EN-SVM model was not able to select a number of features that was under the predefined value $B$. Results indicate that, while in some cases such as for the Breast Cancer P., Diabetes, Ionosphere, Colorectal, and DLBCL datasets, selecting more features led EN-SVM to have better or equal ACC values on the test set than FS-SVM, in all other cases where EN-SVM selected more than $B$ features, it did not lead to improved ACC values in the test set. Thus our model performed better even with fewer features being selected.}

\begin{table}[h]
    \centering
    \renewcommand{\arraystretch}{1.2} 
    \begin{tabular}{lccc|cccc}
        & \multicolumn{3}{c|}{FS-SVM} & \multicolumn{4}{c}{EN-SVM} \\
        \hline
        \hline
        Dataset & ACC Train & ACC Test & Feat. sel. & $\alpha$ & ACC Train & ACC Test & Feat. sel. \\
        \hline
        Breast Cancer D. & \textbf{98.4} & \textbf{93.7} & 5 & $2^9$ & 96.2 & 93.0 & 5 \\
        Breast Cancer P. & \textbf{84.8} & 77.6 & 5 & $2^5$ & 80.7 & \textbf{79.6} & \underline{16} \\
        Breast Cancer W. & \textbf{97.1} & \textbf{96.5} & 5 & $2^{10}$ & 95.5 & 95.3 & \underline{6} \\
        Cleveland & \textbf{85.1} & \textbf{78.7} & 5 & $2^8$ & 78.8 & 70.7 & \underline{8} \\
        Diabetes & \textbf{77.6} & \textbf{76.6} & 5 & $2^9$ & 76.9 & \textbf{76.6} & 4 \\
        German & \textbf{75.3} & \textbf{73.2} & 5 & $2^8$ & 74.4 & 71.6 & \underline{12} \\
        Ionosphere & \textbf{89.0} & \textbf{89.8} & 5 & $2^8$ & 88.2 & \textbf{89.8} & \underline{7} \\
        Parkinsons & \textbf{89.0} & \textbf{85.7} & 5 & $2^7$ & 86.3 & 83.7 & 5 \\
        Wholesale & \textbf{92.1} & \textbf{90.0} & 5 & $2^7$ & 91.5 & 89.1 & 5 \\
        \hline
    \end{tabular}
    \caption{Strict cardinality control (FS) vs.~elastic net penalty control (EN) on the first (small) datasets}
    \label{tab:E-NET1}
\end{table}
\begin{table}[h]
    \centering
    \renewcommand{\arraystretch}{1.2} 
    \begin{tabular}{lccc|cccc}
        & \multicolumn{3}{c|}{FS-SVM} & \multicolumn{4}{c}{EN-SVM} \\
        \hline
        \hline
        Dataset & ACC Train & ACC Test & Feat. sel. & $\alpha$ & ACC Train & ACC Test & Feat. sel. \\
        \hline
        Arrhythmia & \textbf{80.8} & \textbf{75.2} & 10 & $2^9$ & 72.6 & 74.3 & 10 \\
        Colorectal & \textbf{100.0} & \textbf{75.0} & 10 & $2^7$ & 86.9 & \textbf{75.0} & \underline{11} \\
        DLBCL & \textbf{100.0} & 95.0 & 10 & $2^7$ & 98.2 & \textbf{100} & \underline{27} \\
        Lymphoma & \textbf{100.0} & \textbf{95.8} & 10 & $2^8$ & 93.6 & 83.3 & \underline{14} \\
        Madelon & {62.6} & \textbf{63.1} & 10 & $2^9$ & \textbf{65.4} & 61.2 & \underline{84} \\
        Mfeat & \textbf{100} & \textbf{100} & 10 & $2^{10}$ & 99.8 & 99.8 & \underline{17} \\
        \hline
    \end{tabular}
    \caption{Strict cardinality control (FS) vs.~elastic net penalty control (EN) on the second (larger) datasets}
    \label{tab:E-NET2}
\end{table}


\section{Conclusions}\label{Sec6}

In this study, we addressed the NP-hard feature selection problem in linear SVMs under a cardinality constraint, which ensures that the model remains interpretable by selecting a limited number of features. Our approach involves formulating the problem as a MIQP problem and introducing novel SDP relaxations to handle its complexity and enhance scalability.

We developed two MIQP formulations: one employing a big-M method and another one using a complementarity constraint. These formulations accommodate the cardinality constraint by integrating binary variables to select features. To solve these challenging formulations, we proposed several tractable SDP relaxations and a decomposed SDP approach, exploiting the sparsity pattern inherent to the problems. This decomposition significantly reduces computational complexity, making the conic relaxations scalable even for datasets with a large number of features.

For practical implementation, we designed both heuristic and exact algorithms based on the SDP relaxations. The heuristic algorithms, informed by the solutions of the relaxations, efficiently search for good upper bounds and feasible solutions. Meanwhile, the exact algorithm iteratively refines these bounds by solving a sequence of MISOCPs, eventually converging to the global optimum.

Numerical experiments on various benchmark datasets demonstrate the effectiveness of our approach. The proposed methods not only outperform \texttt{Gurobi} but also provide competitive classification accuracy. By allowing a flexible adjustment of the number of features, our approach enhances interpretability without significantly compromising predictive performance.

Future work may explore extensions of this framework to non-linear SVMs or Support Vector Regression problems. Moreover, the reformulations and algorithms proposed in this paper can be extended to other linear-SVM-based models where feature selection plays a crucial role in model interpretation and performance, e.g. optimal classification trees with margins~\cite{MARGOT}.

\section*{Acknowledgements}
 Laura Palagi acknowledges financial support from Progetto di Ricerca Medio Sapienza Uniroma1 (2022) - n. RM1221816BAE8A79. 
 Research of Bo Peng supported by the doctoral programme Vienna Graduate School on Computational Optimization, FWF (Austrian Science Fund), Project W1260-N35. \red{All authors are indebted to excellent reviewers' reports who provided significant remarks for improving presentation;  based upon their suggestion, some new experiments have been added as well, corroborating our approach.}

\bibliography{biblio} 

\begin{thebibliography}{}

\bibitem[Agler et~al., 1988]{agler1988positive}
Agler, J., Helton, W., McCullough, S., and Rodman, L. (1988).
\newblock Positive semidefinite matrices with a given sparsity pattern.
\newblock {\em Linear Algebra and its Applications}, 107:101--149.

\bibitem[Agor and {\"O}zalt{\i}n, 2019]{agor2019feature}
Agor, J. and {\"O}zalt{\i}n, O.~Y. (2019).
\newblock Feature selection for classification models via bilevel optimization.
\newblock {\em Computers \& Operations Research}, 106:156--168.

\bibitem[Alon et~al., 1999]{alon1999broad}
Alon, U., Barkai, N., Notterman, D.~A., Gish, K., Ybarra, S., Mack, D., and Levine, A.~J. (1999).
\newblock Broad patterns of gene expression revealed by clustering analysis of tumor and normal colon tissues probed by oligonucleotide arrays.
\newblock {\em Proceedings of the National Academy of Sciences}, 96(12):6745--6750.

\bibitem[Angelelli et~al., 2010]{angelelli2010kernel}
Angelelli, E., Mansini, R., and Speranza, M.~G. (2010).
\newblock Kernel search: A general heuristic for the multi-dimensional knapsack problem.
\newblock {\em Computers \& Operations Research}, 37(11):2017--2026.

\bibitem[Angelelli et~al., 2012]{angelelli2012kernel}
Angelelli, E., Mansini, R., and Speranza, M.~G. (2012).
\newblock Kernel search: A new heuristic framework for portfolio selection.
\newblock {\em Computational Optimization and Applications}, 51:345--361.

\bibitem[Aytug, 2015]{aytug2015feature}
Aytug, H. (2015).
\newblock Feature selection for support vector machines using generalized {B}enders decomposition.
\newblock {\em European Journal of Operational Research}, 244(1):210--218.

\bibitem[Berman and Shaked-Monderer, 2003]{berman_completely_2003}
Berman, A. and Shaked-Monderer, N. (2003).
\newblock {\em Completely Positive Matrices}.
\newblock World Scientific.

\bibitem[Bertsimas and Shioda, 2009]{bertsimas2009algorithm}
Bertsimas, D. and Shioda, R. (2009).
\newblock Algorithm for cardinality-constrained quadratic optimization.
\newblock {\em Computational Optimization and Applications}, 43(1):1--22.

\bibitem[Bertsimas and Van~Parys, 2020]{bertsimas2020sparse}
Bertsimas, D. and Van~Parys, B. (2020).
\newblock Sparse high-dimensional regression.
\newblock {\em The Annals of Statistics}, 48(1):300--323.

\bibitem[Bradley and Mangasarian, 1998]{Bradley1998FeatureSV}
Bradley, P.~S. and Mangasarian, O.~L. (1998).
\newblock Feature selection via concave minimization and support vector machines.
\newblock In {\em International Conference on Machine Learning}.

\bibitem[Chan et~al., 2007]{chan2007direct}
Chan, A.~B., Vasconcelos, N., and Lanckriet, G.~R. (2007).
\newblock Direct convex relaxations of sparse {SVM}.
\newblock In {\em Proceedings of the 24th International Conference on Machine Learning}, pages 145--153.

\bibitem[Cheramin et~al., 2022]{cheramin2022computationally}
Cheramin, M., Cheng, J., Jiang, R., and Pan, K. (2022).
\newblock Computationally efficient approximations for distributionally robust optimization under moment and {W}asserstein ambiguity.
\newblock {\em INFORMS Journal on Computing}, 34(3):1768--1794.

\bibitem[Cortes and Vapnik, 1995]{VapnikSVM1995}
Cortes, C. and Vapnik, V. (1995).
\newblock Support-{v}ector {n}etworks.
\newblock {\em Machine Learning}, 20(3):273--–297.

\bibitem[Dua and Graff, 2017]{UCI2019}
Dua, D. and Graff, C. (2017).
\newblock {UCI} {M}achine {L}earning {R}epository.
\newblock \url{http://archive.ics.uci.edu/ml}.

\bibitem[D’Onofrio et~al., 2024]{MARGOT}
D’Onofrio, F., Grani, G., Monaci, M., and Palagi, L. (2024).
\newblock Margin optimal classification trees.
\newblock {\em Computers \& Operations Research}, 161:106441.

\bibitem[Fan et~al., 2008]{fan2008liblinear}
Fan, R.-E., Chang, K.-W., Hsieh, C.-J., Wang, X.-R., and Lin, C.-J. (2008).
\newblock Liblinear: A library for large linear classification.
\newblock {\em Journal of Machine Learning Research}, 9:1871--1874.

\bibitem[Fung and Mangasarian, 2004]{FungFSNewton}
Fung, G. and Mangasarian, O.~L. (2004).
\newblock A feature selection {N}ewton method for support vector machine classification.
\newblock {\em Computational Optimization and Applications}, 28:185--202.

\bibitem[Gao and Li, 2013]{gao2013optimal}
Gao, J. and Li, D. (2013).
\newblock Optimal cardinality constrained portfolio selection.
\newblock {\em Operations research}, 61(3):745--761.

\bibitem[Ghaddar and Naoum-Sawaya, 2018]{ghaddar2018high}
Ghaddar, B. and Naoum-Sawaya, J. (2018).
\newblock High dimensional data classification and feature selection using support vector machines.
\newblock {\em European Journal of Operational Research}, 265(3):993--1004.

\bibitem[Guastaroba and Speranza, 2012]{guastaroba2012kernel}
Guastaroba, G. and Speranza, M.~G. (2012).
\newblock Kernel search for the capacitated facility location problem.
\newblock {\em Journal of Heuristics}, 18(6):877--917.

\bibitem[{Gurobi Optimization, LLC}, 2024]{gurobi}
{Gurobi Optimization, LLC} (2024).
\newblock {Gurobi Optimizer Reference Manual}.

\bibitem[Guyon et~al., 2002]{guyon2002gene}
Guyon, I., Weston, J., Barnhill, S., and Vapnik, V. (2002).
\newblock Gene selection for cancer classification using support vector machines.
\newblock {\em Machine Learning}, 46:389--422.

\bibitem[Hastie et~al., 2004]{hastie2004entire}
Hastie, T., Rosset, S., Tibshirani, R., and Zhu, J. (2004).
\newblock The entire regularization path for the support vector machine.
\newblock {\em Journal of Machine Learning Research}, 5(Oct):1391--1415.

\bibitem[Hastie et~al., 2009]{hastie2009elements}
Hastie, T., Tibshirani, R., and Friedman, J.~H. (2009).
\newblock {\em The Elements of Statistical Learning: Data Mining, Inference, and Prediction}.
\newblock Springer, 2nd edition.

\bibitem[Kamath and Liu, 2021]{kamath2021explainable}
Kamath, U. and Liu, J. (2021).
\newblock {\em Explainable Artificial Intelligence: An Introduction to Interpretable Machine Learning}.
\newblock Springer.

\bibitem[Kanzow et~al., 2021]{kanzow2021augmented}
Kanzow, C., Raharja, A.~B., and Schwartz, A. (2021).
\newblock An augmented lagrangian method for cardinality-constrained optimization problems.
\newblock {\em Journal of Optimization Theory and Applications}, 189(3):793--813.

\bibitem[Labb{\'e} et~al., 2019]{labbe2019mixed}
Labb{\'e}, M., Mart{\'\i}nez-Merino, L.~I., and Rodr{\'\i}guez-Ch{\'\i}a, A.~M. (2019).
\newblock Mixed integer linear programming for feature selection in support vector machine.
\newblock {\em Discrete Applied Mathematics}, 261:276--304.

\bibitem[Lee et~al., 2020]{lee2020mixed}
Lee, I.~G., Zhang, Q., Yoon, S.~W., and Won, D. (2020).
\newblock A mixed integer linear programming support vector machine for cost-effective feature selection.
\newblock {\em Knowledge-based systems}, 203:106145.

\bibitem[Maldonado et~al., 2014]{maldonado2014feature}
Maldonado, S., P{\'e}rez, J., Weber, R., and Labb{\'e}, M. (2014).
\newblock Feature selection for support vector machines via mixed integer linear programming.
\newblock {\em Information Sciences}, 279:163--175.

\bibitem[Nguyen and De~la Torre, 2010]{nguyen2010optimal}
Nguyen, M.~H. and De~la Torre, F. (2010).
\newblock Optimal feature selection for support vector machines.
\newblock {\em Pattern Recognition}, 43(3):584--591.

\bibitem[Rudin et~al., 2022]{rudin2022interpretable}
Rudin, C., Chen, C., Chen, Z., Huang, H., Semenova, L., and Zhong, C. (2022).
\newblock Interpretable machine learning: Fundamental principles and 10 grand challenges.
\newblock {\em Statistic Surveys}, 16:1--85.

\bibitem[Shipp et~al., 2002]{shipp2002diffuse}
Shipp, M.~A., Ross, K.~N., Tamayo, P., Weng, A.~P., Kutok, J.~L., Aguiar, R.~C., Gaasenbeek, M., Angelo, M., Reich, M., Pinkus, G.~S., et~al. (2002).
\newblock Diffuse large b-cell lymphoma outcome prediction by gene-expression profiling and supervised machine learning.
\newblock {\em Nature Medicine}, 8(1):68--74.

\bibitem[Tillmann et~al., 2024]{tillmann2024cardinality}
Tillmann, A.~M., Bienstock, D., Lodi, A., and Schwartz, A. (2024).
\newblock Cardinality minimization, constraints, and regularization: a survey.
\newblock {\em SIAM Review}, 66(3):403--477.

\bibitem[Vapnik, 1999]{vapnik1999nature}
Vapnik, V. (1999).
\newblock {\em The Nature of Statistical Learning Theory}.
\newblock Springer Science \& Business Media.

\bibitem[Wang et~al., 2006]{wang2006doubly}
Wang, L., Zhu, J., and Zou, H. (2006).
\newblock The doubly regularized support vector machine.
\newblock {\em Statistica Sinica}, 16(2):589--615.

\bibitem[Weston et~al., 2003]{weston2003use}
Weston, J., Elisseeff, A., Sch{\"o}lkopf, B., and Tipping, M. (2003).
\newblock Use of the zero norm with linear models and kernel methods.
\newblock {\em Journal of Machine Learning Research}, 3:1439--1461.

\bibitem[Yuan et~al., 2010]{yuan2010comparison}
Yuan, G.-X., Chang, K.-W., Hsieh, C.-J., and Lin, C.-J. (2010).
\newblock A comparison of optimization methods and software for large-scale l1-regularized linear classification.
\newblock {\em Journal of Machine Learning Research}, 11:3183--3234.

\bibitem[Zou and Hastie, 2005]{zou2005regularization}
Zou, H. and Hastie, T. (2005).
\newblock Regularization and variable selection via the elastic net.
\newblock {\em Journal of the Royal Statistical Society Series B}, 67(2):301--320.

\end{thebibliography}
\end{document}